\documentclass[a4paper,11pt]{amsart}
\usepackage[english]{babel}
\usepackage{setspace}
\usepackage[a4paper,margin=1in]{geometry}
\usepackage{amsmath}
\usepackage{bm}
\usepackage{amssymb}
\usepackage{amsthm}
\usepackage{graphicx}
\usepackage{verbatim}
\numberwithin{equation}{section}
\usepackage{float}
\usepackage[section]{placeins}
\newtheorem{theorem}{Theorem}[section]
\newtheorem{lemma}{Lemma}[section]
\newtheorem{pro}{Proposition}[section]
\newtheorem{definition}{Definition}[section]
\newtheorem{remark}{Remark}[section]
\newtheorem{assume}{Assumption}[section]

\usepackage{hyperref}
\hypersetup{hidelinks} 
\usepackage{dsfont}
\allowdisplaybreaks[1]
\usepackage{enumerate}
\usepackage{color}
\usepackage{booktabs}
\usepackage{authblk}
\usepackage{amsaddr}
\newcommand{\blue}[1]{{\color{black}#1}}

\usepackage{natbib}
 \bibpunct[, ]{[}{]}{,}{n}{}{,}%

\title[Retirement under Age-Dependent Mortality]{Optimal retirement choice under age-dependent force of mortality}

\author[1]{ Giorgio Ferrari$^\dagger$ \quad and  \quad Shihao Zhu$^\ddagger$}
\thanks{$^\dagger$Center for Mathematical Economics (IMW), Bielefeld University, Universit{\"a}tsstrasse 25, 33615, Bielefeld, Germany (giorgio.ferrari@uni-bielefeld.de).}
\thanks{$^\ddagger$Institute of Insurance Science, Ulm University, Helmholtzstr.\ 20, 89069 Ulm, Germany (shihao.zhu@uni-ulm.de).}

\date{\today}

\begin{document}
%\tableofcontents
\newpage
\maketitle

\begin{abstract}

This paper examines {optimal investment, consumption, and retirement timing} under an age-dependent force of mortality. We formulate the optimization problem as a combined stochastic control and optimal stopping problem with a random time horizon, featuring wealth, labor income, and the force of mortality. To address this problem, we transform it into its dual form, which is a finite time horizon, three-dimensional optimal stopping problem with joint dynamics. We establish the existence of an optimal retirement boundary that separates the state space into continuation and stopping regions. {We derive regularity of the optimal stopping value function, prove that the dual boundary is locally Lipschitz, and characterize the boundary as the unique solution to a nonlinear integral equation, which is solved numerically.} In the original coordinates, the agent retires whenever her wealth exceeds an age-, labor-income-, and force-of-mortality-dependent transformation of the optimal stopping boundary. We also provide numerical illustrations of the optimal strategies, including sensitivity analyses of the optimal retirement boundary with respect to the relevant model parameters.
\end{abstract}

\vspace{2mm}


\vspace{2mm}

\vspace{2mm}

\vspace{2mm}

\section{Introduction}

The timing of retirement is a crucial financial decision that individuals must face as they approach later stages of life. However, determining the appropriate moment to retire is a complex undertaking, as it is influenced by numerous factors. One of the most prominent and evident predictors of retirement decisions is an individual's chronological age, as highlighted in empirical studies such as \citet{wang2014psychological}. Additionally, research indicates that individuals typically postpone retirement until they possess the necessary financial means to do so, with retirement timing being linked to both wealth and labor income, as explored in \citet{honig1998married}. In line with rational decision-making principles and the life-cycle model in economics, subjective life expectancy ({how long one expects to live}) also plays a significant role in retirement timing. \citet{szinovacz2014recession} used data from the Health and Retirement Study and found that expected mortality risk influenced retirement plans, and the effects were especially strong for expectations to work beyond age 65. Furthermore, factors like health status and {family} caregiving responsibilities contribute to the complexity of retirement decision-making. {For a comprehensive examination of these factors, see the review by \citet{fisher2016retirement}.}

In addition to the aforementioned empirical studies, the field of optimal retirement {decision-making} has been enriched by various theoretical models from a financial perspective. The historical origins of the optimal retirement time problem can be traced back to seminal works such as \citet{jin2006disutility}, \citet{farhi2007saving}, \citet{choi2008optimal}, \citet{dybvig2010lifetime}, and \citet{dybvig2011verification}. These foundational contributions paved the way for investigating the optimal investment and consumption behavior of individuals facing retirement decisions, resulting in significant advancements across various contexts. Notable advancements include the introduction of features like mandatory retirement dates and early retirement options (\citet{yang2018optimal}, \citet{jeon2023horizon}), the consideration of consumption ratcheting (\citet{jeon2020optimal}), the exploration of partial information (\citet{chen2021optimal}), the incorporation of habit persistence (\citet{chen2018optimal, he2022optimal, guan2025retirement}), the examination of return ambiguity in risky asset prices (\citet{park2023robust}), the analysis of non-Markovian environments (\citet{yang2021optimal}), and the study of heterogeneous consumption patterns involving basic and luxury goods (\citet{jang2024optimal}).

However, an important characteristic shared by much of the existing literature is the absence of consideration for mortality risk or the assumption of a constant force of mortality. This requirement, although common, can impose significant restrictions on the modeling framework.  In reality, mortality has long been studied in the fields of mathematics and demography. Abraham De Moivre, in 1725, proposed perhaps the earliest mathematical mortality model, suggesting that the probability of surviving from birth to age $x$ follows a linear function of age. Moreover, the most successful and influential mortality law is known as Gompertz's law, formulated by \citet{gompertz1825nature}. This law describes a function with exponentially increasing force of mortality and has found application in insurance economics as well (see, e.g., \citet{milevsky2007annuitization} and \citet{de2019free} for the optimal annuitization problem with a Gompertz force of mortality).  %\textcolor{red}{comparison with  \citet{huang2012optimal}, \citet{bloom2014optimal}, \citet{mao2014optimal}  }

In this paper, we aim to quantitatively examine the interplay between age, wealth, labor income, and mortality on the optimal retirement time. We {incorporate} an age-dependent force of mortality, specifically the Gompertz model, to comprehensively analyze these influences. Our study focuses on an individual who confronts a mandatory retirement date but possesses the option to retire earlier, assuming irreversibility in retirement decisions. The individual can allocate her wealth to consumption and invest in a risky asset. Additionally, the labor income of the individual is stochastic until retirement, and post-retirement utility increases due to the availability of more leisure time.

As highlighted by \citet{chen2021retirement}, ``the optimal retirement problem of a sophisticated individual with an age-dependent force of mortality has not been solved analytically.'' Consequently, the primary objectives of our paper are twofold: first, to determine optimal consumption and portfolio-choice strategies considering realistic mortality risks stemming from an age-dependent force of mortality, and second, to fully characterize the optimal retirement time. By addressing these objectives, our study intends to fill the existing gap in the literature. {More specifically, \citet{chen2021retirement} solve the age-dependent-mortality cases for precommitted and naive agents but obtain the sophisticated solution only under constant mortality, whereas our stopping-time formulation retains Gompertz mortality. \citet{jang2024optimal} allow general time-dependent discounting but study nonhomothetic utility over basic and luxury consumption; we instead use homothetic power utility and provide a complete verification and free-boundary characterization.}

Notably, our findings, which are collected in Section \ref{sec5} and Section \ref{sec6}, are somewhat intuitively convincing. For example, we show that it is optimal to retire when the agent's wealth first reaches an endogenously determined boundary surface, which depends on the agent's age, labor income and force of mortality. Intuitively, if the agent is sufficiently rich (her wealth exceeds the corresponding boundary), then retirement should be performed immediately; otherwise, the agent continues working until wealth reaches the retirement boundary. {In the numerical example with $\gamma>1$, the consumption--wealth ratio and the risky-asset share fall discretely at retirement. After retirement, the risky-asset share is constant and equal to the Merton share, whereas the consumption--wealth ratio depends deterministically on the current force of mortality.}
 We also provide some interesting economic implications of the optimal retirement boundary through a numerical study in Section \ref{sec6}.
%\vspace{-0.2mm}

\vspace{-0.5\baselineskip}
\subsection{Overview of the mathematical analysis}
{From a mathematical point of view, we model the previous problem as a random-time-horizon stochastic control problem with discretionary stopping. Its stochastic state variables are wealth $W$ and labor income $Y$, while the force of mortality $M$ follows a deterministic Gompertz trajectory.} The dominant feature of the wealth process $W$ is that it is not the same before and after the retirement time $\tau$. The agent's aim is to choose consumption rate $c$, portfolio $\pi$, and the retirement time $\tau$ in order to maximize the total expected utility of consumption $c$, up to the random death time $\eta$.

In the literature, the derivation of the Hamilton-Jacobi-Bellman (HJB) equation is a frequently employed approach when analyzing {a} stochastic optimization problem with Markovian processes. %In our specific case, the corresponding HJB equation turns out to be an involved combination of an HJB-type equation (reflecting the investment and consumption optimization) with a variational inequality (reflecting the retirement optimization).
However, studying the properties of the value function and the optimal strategies, particularly the optimal stopping boundary, proves to be challenging when directly studying the HJB equation. {This difficulty arises from the interaction of controlled wealth, stochastic labor income, the retirement option, and time-varying mortality.} In order to tame the intricate mathematical structure of our problem, where the consumption and portfolio choices nontrivially interact with the retirement decision, we combine a duality and a free-boundary approach, and proceed in our analysis as follows.

\textbf{Step 1.} {First, we conduct successive transformations (see Section \ref{sec3}) that connect the original stochastic control-stopping problem (with value function $V$) with its dual problem (with value function $J$) by martingale and duality methods similar to \citet{karatzas2000utility}. The dual problem is written in the variables $(Z,M,Y)$: $M$ is deterministic, the drift of $Z$ depends on $M$, and $Z$ and $Y$ share the same Brownian motion. Homogeneity of power utility and a change of measure reduce $(Z,M,Y)$ to $(X,M)$, with value function $\widetilde J$.}

\textbf{Step 2.} {We analyze the reduced stopping problem in five stages. First, we subtract the payoff from immediate retirement and work with the value function $\widehat J(t,x,m)$. Second, the augmented representation in $(t,x,m)$ embeds all Gompertz trajectories in one time-homogeneous Markov family; it is used for mortality comparisons and for describing the boundary surface $(t,m)\mapsto b(t,m)$. Third, after fixing one Gompertz characteristic, mortality is absorbed into calendar time and the restriction, denoted by $\mathcal J(t,x)$, becomes a one-dimensional time-inhomogeneous stopping problem. Fourth, we establish the value-function and boundary regularity needed for the free-boundary analysis. Finally, these results yield the nonlinear integral characterization used in the primal recovery and numerical analysis. Thus, the logical flow is the original control--stopping problem $\Longrightarrow$ the dual stopping problem $\Longrightarrow$ the fixed-characteristic stopping problem $\Longrightarrow$ the boundary characterization $\Longrightarrow$ the primal policies.}

%It is also worth pointing out that the force of mortality process $M$ is a deterministic function of time $t$, which helps us study the regularity of $\widehat{J}$ (given by the difference of $\widetilde{J}$ and the smooth payoff of immediate stopping; see (\ref{4-3})).
%\blue{The force of mortality $M$ is deterministic. We retain its current level as a state parameter in order to describe, within one time-homogeneous Markov family, agents with different current mortality rates and to obtain comparative statics with respect to mortality. Along any fixed Gompertz trajectory, however, $M$ can be absorbed into calendar time and the genuinely stochastic state is one-dimensional. Thus, the formulation in $(t,X,M)$ is useful but is not intrinsically three-dimensional.}

{The two representations therefore have distinct roles. The augmented formulation $\widehat J(t,x,m)$ supplies comparisons across current mortality levels, whereas the $\mathcal J(t,x)$ formulation supplies one-dimensional parabolic regularity. For the latter step, Appendix \ref{prolipschitz} applies the one-dimensional boundary-regularity result in Theorem 4.3 of \citet{de2019lipschitz} after verifying its assumptions. The exact comparison across Gompertz characteristics and the terminal-boundary argument are specific to the present model.} {In particular, the boundary is locally Lipschitz along each fixed Gompertz characteristic. This regularity yields continuity of the optimal stopping time and supports the nonlinear integral characterization, while the mortality comparison is used to establish the terminal boundary limit.}

%As a matter of fact,

%the process ($X,M$) is a degenerate diffusion process (in the sense that the differential operator of ($X,M$) is a degenerate parabolic operator) so that the study of the regularity of $\widehat{J}$ in the interior of its continuation region cannot hinge on classical  results for parabolic PDEs (see, e.g., \citet{friedman1026variational}).

%Although the new introduced auxiliary optimal stopping problem can

 % In fact, since  we were unable to establish continuity for the mapping $(t,m) \mapsto b(t,m)$, even if monotonicity properties of the latter could be easily proved.

 %However,

%On the other hand, we show that the optimal boundary is in fact a locally Lipschitz-continuous function of time $t$ and force of mortality $m$,  which is a sharp result

\textbf{Step 3}. %After proving the strict convexity of ${J}$,
We come back to the original coordinates' system and via the duality relations we obtain the optimal consumption and portfolio policies, as well as the optimal retirement time, in terms of the optimal stopping boundary and value function (cf. Section \ref{sec5}).

Overall, we believe that the contributions of this paper are the following. From a mathematical point of view, we provide a rigorous theoretical analysis of the optimal retirement time in terms of the optimal boundary $\widehat{b}$ (see (\ref{bhat})). {More specifically, we provide a complete analytical characterization of the value function and the optimal strategies for the present mortality-dependent retirement problem.} Furthermore, we believe that the dual optimal stopping problem, studied as a device to characterize the optimal solution of the optimal retirement time problem, is of interest {in its own right}. By performing a thorough analysis {of} the regularity of $\mathcal{J}$ (an equivalent version of the dual value function) and of the free boundary, we are able to provide a complete characterization of the optimal retirement time strategy through a nonlinear integral equation. The analysis in this paper is completed by solving numerically the integral equation and studying its sensitivity to variations in the model's parameters, thus extending the results of the related optimal retirement time model.
\vspace{-0.5\baselineskip}

\subsection{Organization of the paper}
The rest of the paper is organized as follows. In Section \ref{sec2}, we introduce the optimal retirement time decision model with an age-dependent force of mortality. We transform the original stochastic control-stopping problem into a pure stopping problem in Section \ref{sec3}, while in Section \ref{sec4} we study the optimal stopping problem. In Section \ref{sec5}, we provide the optimal retirement time boundary, {optimal consumption and portfolio policies} in primal variables, and in Section \ref{sec6} we present a numerical study and provide some economic implications. {Section \ref{sec7} concludes and discusses the possible extensions.} Appendix A collects the proofs of some results of Sections \ref{sec3}, \ref{sec4} and \ref{sec5}, whereas Appendix B gives the details of the numerical method used in Section \ref{sec6}.

\section{Setting and problem formulation}\label{sec2}
\subsection{The age-dependent force of mortality}

\leavevmode
{
Fix a mandatory retirement time $T>0$ and an initial time $t\in[0,T]$. For mortality risk, let $({\Omega}^M,{\mathcal{F}^M},{\mathbb{P}^M})$ be a complete probability space supporting an exponential random variable $\mathcal{E}$, so that $\mathbb{P}^M(\mathcal{E}>v)=e^{-v}$ for $v\geq0$. 
}

 Consider an agent whose force of mortality evolves according to the standard Gompertz model (see, e.g., \citet{gompertz1825nature}). The force of mortality process $M:=\{M_s, s\geq t\}$ thus follows the dynamics
 \begin{align}\label{2-1}
d M_s = aM_s ds ,  \ M_t =m>0,
\end{align}
with $a>0$. Notice that the usual form for the Gompertz model is $M_s=a_0^{-1}e^{(x_0+s-t-m_0)/{a_0}}$, hence we are using $a=\frac{1}{a_0}$ and $m=a_0^{-1}e^{(x_0-m_0)/{a_0}}$. Here, $x_0$ denotes the agent's age at initial time $t$, $m_0$ is called the modal value, and $a_0$ is the dispersion parameter of the Gompertz model. This parsimonious specification has a long history of calibration to population mortality data.

{
Following the Cox-process construction in \citet[Chapter 2.3.1]{aksamit2017enlargement}, conditional on the agent being alive at $t$, we define the remaining death time by
\begin{align*}
\eta:=\inf\bigg\{s\geq t:\int_t^s M_u du\geq\mathcal{E}\bigg\}.
\end{align*}
For any $s\geq t$, this construction gives the survival probability (cf.\ \citet[Chapter 2.3.2, Lemma 2.25, p. 42]{aksamit2017enlargement})
\begin{align*}
\mathbb{P}^M(\eta>s)=\mathbb{P}^M\bigg(\mathcal{E}>\int_t^sM_u du\bigg)=\exp\bigg\{-\int_t^sM_u du\bigg\}.
\end{align*} 
}

%\begin{sloppypar}

\subsection{The financial market and labor income}
%\leavevmode
{
Let $({\Omega},{\mathcal{F}},{\mathbb{P}})$ be a complete probability space endowed with the augmented filtration $\mathbb{F}:=\{\mathcal{F}_s,s\geq0\}$ generated by a standard Brownian motion $B:=\{B_s,s\geq0\}$. We work on the completion of the product probability space
\begin{align*}
(\Omega^M\times\Omega,\mathcal{F}^M\otimes\mathcal{F},\mathbb{P}^M\otimes\mathbb{P}).
\end{align*}
We denote this completed product space by $\mathcal{M}$ and extend the mortality and financial random variables canonically to it. In particular, $\mathcal{E}$, and hence $\eta$, is independent of $\mathcal{F}_\infty$. We denote expectation under the completed product measure by $\mathbb{E}^{\mathcal{M}}$.
}
We assume that the agent invests in a financial market with two assets. One of them is a risk-free bond, whose price $S^0:=\{ S^0_s, s \geq t \}$ evolves as 
\begin{align*}
dS^0_s=rS^0_s ds, \quad S^0_t=s^0>0,
\end{align*}
where $r>0$ is a constant risk-free rate. The second one is a stock, whose price is denoted by $S:=\{S_s, s\geq t\}$ and it satisfies the stochastic differential equation
\begin{align*}
dS_s &= \mu S_s ds + \sigma S_s dB_s, \quad S_t=s>0,
\end{align*}
where $\mu \in \mathbb{R}$ and $\sigma>0$ are given constants.
\enlargethispage{3pt}

Let $\tau\in[t,T]$ be an $\mathbb{F}$-stopping time representing the time at which the agent chooses to retire. The potential labor income process $Y:=\{Y_s,t\leq s\leq T\}$ is spanned by the market and follows a geometric Brownian motion, similar to e.g.\ \citet{bodie2004optimal} and \citet{dybvig2010lifetime},\footnote{Since most salaries are not directly tied to financial market risk, assuming perfectly hedgeable labor income is not fully realistic. We nevertheless adopt this simplification to focus on the impact of mortality risk on retirement, following the tractable framework of \citep{bodie2004optimal, dybvig2010lifetime}. Extending the model to partially unhedgeable income risk (e.g., with an additional Brownian motion) is left for future research.}
\begin{align}\label{2-4-1}
Y_s=y e^{(\mu_y-\frac{1}{2}\sigma_y^2)(s-t)+\sigma_y(B_s-B_t)}, \quad t\leq s\leq T, \quad Y_t=y>0,
\end{align}
where $\mu_y\in\mathbb{R}$ and $\sigma_y>0$ are constants representing the instantaneous growth rate and volatility of potential labor income, respectively. The agent receives income at rate $Y_s$ only while $s<\tau$.

We define the market price of risk $\theta:= \frac{\mu-r}{\sigma}$ and the state-price-density process $\xi_{s,t}:=e^{-r(s-t)-\theta (B_s-B_t) -\frac{1}{2}\theta^2 (s-t)}$. Since the market is complete and the labor income is perfectly correlated with the market, {the present value of future labor income $g_s$ at time $s<T$ (also called ``human capital" as in \citet{dybvig2010lifetime}), under the assumption that the agent does not choose early retirement and the agent is alive before $T$,} is given by
\begin{align}\label{b(s)}
g_s&= \mathbb{E} \bigg[\int^{T}_s \xi_{u,s} Y_u du\bigg| \mathcal{F}_s \bigg] =Y_s \mathbb{E} \bigg[\int^{T}_s \xi_{u,s} \frac{Y_u}{Y_s} du\bigg| \mathcal{F}_s\bigg]=Y_s  \bigg[\int^{T}_s e^{-(r-\mu_y+\sigma_y \theta)(u-s)}  du\bigg] \nonumber\\
&=
\begin{aligned}
&Y_s \frac{1-e^{-\kappa(T-s)}}{\kappa},
%&Y_s (T-s), \ &\text{if} \ \kappa = 0
\end{aligned}
\end{align} 
where $\kappa:=r-\mu_y+\sigma_y\theta$ is the effective discount rate for labor income, which is assumed to be strictly positive (see also \citet{dybvig2011verification} and \citet{guan2025retirement} for a similar requirement). 

 Then, for $s\leq T$, we define 
\begin{align*}
q_s:=\frac{1-e^{-\kappa(T-s)}}{\kappa},
\end{align*}
so that $g_s=q_sY_s$.

The agent also consumes from her wealth, while investing in the financial market. Denoting by $\pi_s$ the amount of wealth invested in the stock at time $s$, the agent then also chooses the rate of spending in consumption $c_s$ at time $s$. Therefore, the agent's wealth $W:=\{W^{c,\pi,\tau}_s, s\geq t\}$ evolves as 
\begin{align}\label{sde}
d W^{c,\pi, \tau}_s =  [\pi_s(\mu-r)+rW^{c,\pi,\tau}_s-c_s + Y_s \blue{\mathds{1}_{\{s<\tau\}}}]ds +\pi_s \sigma dB_s, \quad W_t^{c,\pi,\tau}= w.
\end{align}
In the following, we shall simply write $W$ to denote $W^{c,\pi, \tau}$, where needed.

\subsection{The optimization problem}
 Here and in the sequel, we write \blue{$\mathcal{O}:=[0,T] \times \mathbb{R}^3_+$}, \blue{$\mathcal{O}':= \{(t,w,m,y)\in [0,T] \times \mathbb{R} \times \mathbb{R}^2_+| w>-y q_t \}$ } and $\mathcal{U}:=[0,T] \times \mathbb{R}^2_+$, with $\mathbb{R}^n_+:=(0,\infty)^n$ for each integer $n\geq 1$. We denote by $\mathcal{S}_{t,s}$ the class of $\mathbb{F}$-stopping times $\tau: \Omega \to [t,s]$ for $t \leq s \leq T$, and let $\mathcal{S} := \mathcal{S}_{t,T}$. Then we introduce the set of admissible strategies as follows.

\begin{definition}\label{admissiblecontrol}
Let $(t,w,m,y) \in \mathcal{O}'$ be given and fixed. The triplet of choices $(c,\pi, \tau)$ is called an \textbf{admissible strategy} for $(t,w,m,y)$, and we write $(c,\pi, \tau) \in \mathcal{A}(t,w,m,y)$, if it satisfies the following conditions: 
\begin{enumerate}[(i)]
	\item  $c$ and $\pi$ are progressively measurable with respect to $\mathbb{F}$, $\tau \in \mathcal{S}$;
	\item $c_s \geq 0$ for all $s \geq t$ and {$\int_t^\Theta(c_s+|\pi_s|^2)ds <\infty$ $\mathbb{P}$-a.s., for every finite $\Theta>t$;} 
	\item $W^{c,\pi, \tau}_s \geq   -g_s\mathds{1}_{\{s< \tau\}}$ for all $s\geq t$, where $g_s$ is defined in (\ref{b(s)}) if $s <T$ and $g_s=0$ if $s \geq T$.
\end{enumerate}
\end{definition}
\textcolor{black}{Condition (iii) means that, before retirement, the agent may have negative financial wealth, but its amount cannot exceed the present value of future labor income. After retirement, when labor income is no longer available, the constraint reduces to non-negativity of net financial wealth. }

The preferences of the agent are described through a power utility function. Due to a larger amount of free time when retired, the agent assigns a higher utility value to consumption during retirement, represented by a constant weighting $K>1$. The utility of the agent for consumption $c_s$ at time $s$ is 
\begin{align*}
U(c_s)= \frac{1}{1-\gamma}\Big[(\mathds{1}_{\{s <\tau\}}+\mathds{1}_{\{\tau \leq s\}}K) \cdot c_s\Big]^{1-\gamma}, \ \gamma \neq 1, \ \gamma>0.
\end{align*}
 Specifically, we call $U_1$ the utility function before retirement, i.e., for $s<\tau$,
\begin{align*}
U_1(c_s):=\frac{1}{1-\gamma}c_s^{1-\gamma};
\end{align*}
on the other hand, we call $U_2$ the utility function after retirement, i.e., for $s\geq \tau$,
\begin{align*}
U_2(c_s):=\frac{1}{1-\gamma}(K \cdot c_s)^{1-\gamma}.
\end{align*}
\blue{Throughout the paper, when $\gamma>1$, an expected-utility criterion is allowed to take the value $-\infty$ under the standard extended-real convention. This convention applies to all subsequent value functions and duality identities and will not be repeated.}

{\color{black}
From the perspective of time $t$, conditional on survival to $t$, the agent maximizes over all $(c,\pi,\tau)\in\mathcal{A}(t,w,m,y)$ the expected  utility functional
\begin{align*}
\mathbb{E}^{\mathcal{M}}\bigg[\int_t^{\eta}e^{-\rho(s-t)}U(c_s)ds\bigg|\mathcal{F}_t\bigg],
\end{align*}
where $\rho>0$ is a constant representing the discount rate. Using the Fubini--Tonelli theorem (applied to $U$ when $\gamma<1$ and to $-U$ when $\gamma>1$), together with the independence of $\eta$ and $\mathcal{F}_\infty$, mortality risk can be integrated out as follows:
\begin{align*}
&\mathbb{E}^{\mathcal{M}}\bigg[\int_t^{\eta}e^{-\rho(s-t)}U(c_s)ds\bigg|\mathcal{F}_t\bigg]\\
&=\mathbb{E}\bigg[\int_t^\infty e^{-\rho(s-t)}\mathbb{P}^M(\eta>s)U(c_s)ds\bigg|\mathcal{F}_t\bigg]\\
&=\mathbb{E}\bigg[\int_t^\infty e^{-\int_t^s(\rho+M_u)du}U(c_s)ds\bigg|\mathcal{F}_t\bigg].
\end{align*} 
}

Hence, given the Markovian setting, the agent aims at determining
\begin{align}\label{2-6}
V(t,w,m,y)=\sup_{(c, \pi,\tau) \in \mathcal{A}(t,w,m,y) } \mathbb{E}_{t,w,m,y}\bigg[  \int_t^{ \infty } e^{- \int^s_t  (\rho +M_u)du}U(c_s)  ds  \bigg],
\end{align}
where $ \mathbb{E}_{t,w,m,y}$ \blue{denotes} the expectation under $\mathbb{P}$ conditioned on  $W_t=w, M_t=m$ and $Y_t=y$. In the rest of the paper, we shall focus on (\ref{2-6}).

		\section{From control-stopping to pure stopping}\label{sec3}
		
		\subsection{The static budget constraint}
		We now transform the dynamic budget constraint in (\ref{sde}) into a static budget constraint based on the standard local martingale argument and the \blue{optional} sampling theorem in \citet{karatzas2000utility}. We obtain the following budget constraint: 
		\begin{align}\label{3-3}
\mathbb{E}_{t,w,m,y}\big[{\xi_{ \tau \wedge s ,t}}(W_{\tau \wedge s }+g_{  \tau \wedge s })\big]+\mathbb{E}_{t,w,m,y}\bigg[ \int_t^{ \tau \wedge s} \xi_{u,t} c_u du                         \bigg]  \leq w+g_t, \quad \text{if} \ s\geq t,
		\end{align}
		and		\blue{\begin{align}\label{3-4}
\mathbb{E}_{t,w,m,y}\bigg[\xi_{s\vee\tau,\tau}W_{s\vee\tau}+\int_\tau^{s\vee\tau}\xi_{u,\tau}c_udu\bigg|\mathcal F_\tau\bigg]\leq W_\tau,\quad s\geq t,
		\end{align}	
		where $\xi_{u,\tau}:=\xi_{u,t}/\xi_{\tau,t}$ for $u\geq\tau$. The stopped intervals in (\ref{3-3}) and (\ref{3-4}) distinguish the regimes $s<\tau$ and $s\geq\tau$.}

	\subsection{The agent's optimization problem after retirement}
%	\leavevmode
	{\color{black}
For any initial time $t\in[0,T]$, suppose that retirement has already occurred and that the post-retirement state is $(W_t,M_t)=(w,m)$, where $w\geq0$ and $m>0$. We first specify admissibility for the post-retirement problem.
\begin{definition}\label{postadmissible}
A consumption--portfolio pair $(c,\pi)$ is called post-retirement admissible for $(t,w,m)$, and we write $(c,\pi)\in\widehat{\mathcal A}(t,w,m)$, if it satisfies the following conditions:
\begin{enumerate}[(i)]
\item $c$ and $\pi$ are progressively measurable with respect to $\mathbb F$;
\item $c_s\geq0$ for $s\geq t$ and $\int_t^\Theta (c_s+|\pi_s|^2)ds<\infty$ $\mathbb P$-a.s. for every finite $\Theta >t$;
\item the wealth process governed by
\begin{align*}
dW_s=[\pi_s(\mu-r)+rW_s-c_s]ds+\sigma\pi_sdB_s,\quad W_t=w,
\end{align*}
satisfies $W_s\geq0$ for all $s\geq t$.
\end{enumerate}
\end{definition}
Accordingly, the agent's post-retirement value function is
}
\begin{align}\label{3-4-1}
\widehat{V}(t,w,m):&=\sup_{(c, \pi) \in \blue{\widehat{\mathcal{A}}(t,w,m)}  } \mathbb{E}_{t,w,m}\bigg[  \int_t^{\infty } e^{- \int^s_t  (\rho +M_u)du}U_2(c_s)  ds \bigg],
\end{align}
where $M$ \blue{is} defined in (\ref{2-1}) and $\mathbb{E}_{t,w,m}$ denotes the expectation conditioned on $W_t=w$ and $M_t=m$.

%\begin{theorem}
%	For any $\tau \in \mathcal{S}$ and $\mathcal{F}_\tau$-measurable random variable $D$ such that $D(\omega) \in \mathbb{R}_+$ for almost sure $\omega \in \Omega$ as well as 
	
%	there exists a strategy $(c,\pi) \in \mathcal{A}_t$ such that the solution to SDE () with initial condition $W_\tau=D$ satisfies 
%	\begin{align*}		W_\tau^{c,\pi} \geq 0 \quad \forall s \geq \tau,
%	\end{align*}
%	and 
%	\begin{align*}
%		\mathbb{E}_{t,w,m}\bigg[  \int_\tau^{\infty } e^{- \int^s_t  (\rho +M_u)du}U_2(c_s)  ds \bigg]= \mathbb{E}[e^{- \int^s_t  (\rho +M_u)du}\widehat{V}(\tau, D,M_\tau)].
%			\end{align*}
%			(2) For $V$ defined in (), we have 
%			\begin{align}\label{3-11-1}
%V(t,w,m,y)=\sup_{(c,\pi,\tau) \in \mathcal{A}(t,w,m,y)} \mathbb{E}_{t,w,m,y} \bigg[ \int_t^{\tau}  e^{-\int^s_t (\rho+ M_u) du} U_1(c_s)  ds   +e^{-\int^{\tau}_t (\rho+ M_u) du}  \widehat{V}( W_{\tau},M_{\tau})    \bigg].
%\end{align}
%\end{theorem}

\blue{For any $(c,\pi)\in\widehat{\mathcal A}(t,w,m)$, similar to (\ref{3-4})} we have  
\blue{\begin{align*}
\mathbb E_{t,w,m}\bigg[\int_t^\infty\xi_{s,t}c_sds\bigg]\leq w.
\end{align*}
For any Lagrange multiplier $z>0$,}
\begin{align}\label{3-5}
 &\mathbb{E}_{t,w,m} \bigg[    \int_t^{\infty }  e^{-\int^s_t (\rho+ M_u) du} U_2(c_s)  ds  \bigg] \nonumber \\&\leq
\mathbb{E}_{t,w,m} \bigg[    \int_t^{\infty}  e^{-\int^s_t (\rho+ M_u) du} U_2(c_s)  ds  \bigg]   - z      \mathbb{E}_{t,w,m}\bigg[  \int_t^{\infty} \xi_{s,t}c_s ds                         \bigg]  +z w                 \nonumber             \\
&=   \mathbb{E}_{t,w,m}\bigg[    \int_t^{\infty }  e^{-\int^s_t (\rho+ M_u) du} \Big(U_2(c_s) -z P_sc_s\Big) ds  \bigg]+ zw  \nonumber \\
&\leq  \mathbb{E}_{t,w,m}\bigg[    \int_t^{\infty }  e^{-\int^s_t (\rho+ M_u) du} U_2^\ast(z P_s) ds  \bigg]+ zw,
\end{align}
where
\begin{align}\label{3-6}
P_s:= \xi_{s,t}e^{\int^s_t (\rho+ M_u) du } \quad \text{and} \quad U_2^\ast(z):= \sup_{c \geq 0}[U_2(c)-c z].
\end{align}

Let then $Z_s :=z P_s$. By It\^o's formula, we obtain that the dual variable $Z$ satisfies
\begin{align}\label{3-7}
dZ_s=  (\rho-r+ M_s)Z_sds - \theta Z_s dB_s, \quad Z_t= z,
\end{align}
and we set
\begin{align}\label{3-8}
Q(t,z,m):=  \mathbb{E}_{t,z,m}\bigg[    \int_t^{\infty }  e^{-\int^s_t (\rho+ M_u) du} U_2^\ast(Z_s) ds  \bigg].
\end{align}

\begin{pro}\label{Wregularity}
$Q$ is independent of time $t$, it is finite and is strictly convex with respect to $z$, and \blue{both $Q$ and $Q_z$ belong to $C^{2,1}(\mathbb{R}^2_+)$.} Moreover, $Q$ satisfies
\begin{align}\label{3-6-1}
-\widehat{\mathcal{L}}Q=U_2^\ast, 
\end{align}
where 
\begin{align*}
\widehat{\mathcal{L}}Q:= \frac{1}{2}\theta^2 z^2Q_{zz}+(\rho-r+m)zQ_z  +a m Q_m-(\rho+m)Q.
\end{align*}
\end{pro}

\begin{proof}
The proof is given in Appendix \ref{proW}.
\end{proof}
Given that $Q$ is time-independent, with a slight abuse of notation in the sequel we then simply write $Q(z,m)$. \blue{Throughout the paper, $Q_z$ denotes the derivative with respect to the first argument, even when that argument is denoted by $x$, as will occur later.} It is possible to deduce properties of $\widehat{V}$ through those of $Q$ by the following duality relation.
 
\begin{theorem}\label{dualityrelation}
\leavevmode
{\color{black}
For $t\in[0,T]$, $m>0$, and $w\geq0$, the following dual relations hold:
\begin{align*}
\widehat{V}(t,w,m) = \inf_{z >0} [ Q(z,m)+z w ], \quad Q(z,m) =  \sup_{w\geq0}[\widehat{V}(t,w,m)-zw].
\end{align*}
For $w>0$, let $z^\ast:=z^\ast(w,m)$ be the unique solution of $w+Q_z(z^\ast,m)=0$, and let $Z^\ast$ solve (\ref{3-7}) with $Z_t^\ast=z^\ast$. Then, for $s\geq t$,
\begin{align*}
c_s^\ast&=K^{\frac{1-\gamma}{\gamma}}(Z_s^\ast)^{-\frac{1}{\gamma}},\quad
\pi_s^\ast=\frac{\theta}{\sigma}Z_s^{\ast}Q_{zz}(Z_s^\ast,M_s),\quad
W_s^\ast=-Q_z(Z_s^\ast,M_s).
\end{align*}
The pair $(c^\ast,\pi^\ast)$ belongs to $\widehat{\mathcal A}(t,w,m)$ and is optimal for (\ref{3-4-1}).
}

\end{theorem}

\begin{proof}
The proof is given in Appendix \ref{produal}.
\end{proof}

From Theorem \ref{dualityrelation} we see that $\widehat{V}$ is time-independent, so that in the following, with a slight abuse of notation, we simply write $\widehat{V}(w,m)$. 

\subsection{The dual optimal stopping problem}

\blue{By Theorem \ref{dualityrelation}, the deterministic-state post-retirement value $\widehat V$ is continuous on $(0,\infty)^2$ and Borel measurable on $[0,\infty)\times(0,\infty)$ under the extended-real convention. In particular, $\widehat V(0,m)=0$ when $\gamma<1$, whereas $\widehat V(0,m)=-\infty$ and $\widehat V(w,m)\to-\infty$ as $w\downarrow0$ when $\gamma>1$. Hence $\widehat V(W_\tau,M_\tau)$ is $\mathcal F_\tau$-measurable for every admissible retirement time $\tau$.
 For any $(c,\pi,\tau)\in\mathcal A(t,w,m,y)$, the continuation controls after $\tau$ are admissible for the post-retirement problem starting from the random state $(W_\tau,M_\tau)$. Therefore, conditionally on $\mathcal F_\tau$, their continuation payoff is bounded above by $\widehat V(W_\tau,M_\tau)$. Taking expectations and then the supremum gives the upper bound}
\begin{align}\label{3-11-1}
V(t,w,m,y)\leq \sup_{(c,\pi,\tau) \in \mathcal{A}(t,w,m,y)} \mathbb{E}_{t,w,m,y} \bigg[ \int_t^{\tau}  e^{-\int^s_t (\rho+ M_u) du} U_1(c_s)  ds  \nonumber\\ +e^{-\int^{\tau}_t (\rho+ M_u) du}  \widehat{V}( W_{\tau},M_{\tau})    \bigg].
\end{align}

Now, for any $(t,w,m,y) \in \mathcal{O}'$ and Lagrange multiplier $z>0$, from the budget constraint (\ref{3-3}) and (\ref{3-11-1}), recalling $P_s$ in (\ref{3-6}), we have
\begin{align}\label{3-12}
V(t,w,m,y)
\leq&\sup_{(c,\pi,\tau)\in\mathcal A(t,w,m,y)}\mathbb E_{t,w,m,y}\bigg[\int_t^\tau e^{-\int_t^s(\rho+M_u)du}U_1(c_s)ds
+e^{-\int_t^\tau(\rho+M_u)du}\widehat V(W_\tau,M_\tau)\nonumber\\
&-z\xi_{\tau,t}(W_\tau+g_\tau)-z\int_t^\tau\xi_{s,t}c_sds\bigg]+z(w+g_t)\nonumber\\
=&\sup_{(c,\pi,\tau)\in\mathcal A(t,w,m,y)}\mathbb E_{t,w,m,y}\bigg[\int_t^\tau e^{-\int_t^s(\rho+M_u)du}\big(U_1(c_s)-zP_sc_s\big)ds\nonumber\\
&+e^{-\int_t^\tau(\rho+M_u)du}\big(\widehat V(W_\tau,M_\tau)-zP_\tau W_\tau-zP_\tau g_\tau\big)\bigg]+z(w+g_t)\nonumber\\
\leq&\sup_{\tau\in\mathcal S}\mathbb E_{t,z,m,y}\bigg[\int_t^\tau e^{-\int_t^s(\rho+M_u)du}U_1^\ast(Z_s)ds\nonumber\\
&+e^{-\int_t^\tau(\rho+M_u)du}\big(Q(Z_\tau,M_\tau)-Z_\tau g_\tau\big)\bigg]+z(w+g_t),
\end{align}
where $Z_s$ is defined in (\ref{3-7}) and 
\begin{align}\label{3-12-1}
{U}_1^\ast(z):= \sup_{c \geq 0}[U_1(c)-c z].
\end{align}

Hence, defining the (candidate) dual value function
\begin{align}\label{3-13}
J(t,z,m,y):=\sup_{\tau \in\mathcal{S}
} \mathbb{E}_{t,z,m,y}\bigg[    \int_t^{\tau }  e^{-\int^s_t (\rho+ M_u) du}{U}^\ast_1(Z_s) ds \nonumber\\ + e^{-\int^{\tau}_t (\rho+ M_u) du}\Big(Q( Z_{\tau}, M_{\tau})-Z_{\tau}g_\tau \Big)\bigg],
\end{align}
\blue{we obtain a finite-horizon optimal stopping problem for the joint state $(Z,M,Y)$. Here, $M$ has deterministic dynamics, the drift of $Z$ depends on $M$, and $Z$ and $Y$ are driven by the same Brownian motion $B$, as specified in (\ref{3-7}), (\ref{2-1}), and (\ref{2-4-1}).} We shall close the duality gap and show that $J$ is indeed the dual value function in Section \ref{sec5} (see Theorem \ref{verification} below). In Section \ref{sec4}, we perform a detailed probabilistic study of (\ref{3-13}).

\section{Study of the (candidate) dual optimal stopping problem}\label{sec4}
\subsection{Dimensionality reduction}

{The finite-horizon dual problem in (\ref{3-13}) is written in the state variables $(z,m,y)$. Homogeneity of $U$ and a change of measure combine $(z,y)$ into one state $x$.}

\begin{lemma}\label{reduce}
Consider the exponential martingale
\begin{align*}
\zeta := \bigg\{ \zeta_s = \exp\bigg(\int^s_t -\frac{1}{2}(\sigma_y-\theta)^2du+\int^s_t (\sigma_y-\theta)dB_u \bigg):  t\leq s \leq T                          \bigg\}.
\end{align*}
Let $\widehat{\mathbb{P}}$ be the probability measure on $(\Omega, \mathcal{F}_T)$ such that $\frac{d\widehat{\mathbb{P}}}{d \mathbb{P}}= \zeta_T$, and denote by $\widehat{\mathbb{E}}$ the expectation under $\widehat{\mathbb{P}}$. Then we have 
\begin{align}\label{jandj}
	J(t,z,m,y)=zy \widetilde{J}(t,z^{\frac{1}{1-\gamma}}y^{\frac{\gamma}{1-\gamma}},m)\end{align} 
with
\begin{align}\label{4-2}
\widetilde{J}(t,x,m):=\sup_{t \leq \tau \leq T
} \widehat{\mathbb{E}}_{t,x,m}\bigg[    \int_t^{\tau }  e^{-\kappa(s-t)} {U}^\ast_1(X_s) ds+ e^{-\kappa(\tau-t)}  \Big(Q(X_\tau, M_\tau)- q_\tau \Big)\bigg],
\end{align}
and the underlying Markov process $X$ being defined as $ X:=Z^{\frac{1}{1-\gamma}}Y^{\frac{\gamma}{1-\gamma}}$ (with $Z$ and $Y$ as in (\ref{3-7}) and (\ref{2-4-1}), respectively).
\end{lemma}
 \begin{proof}
The proof is given in Appendix \ref{proreduce}.
\end{proof}

We next introduce the dynamic equation of $X$ defined in Lemma \ref{reduce}. By It\^o's formula and Girsanov's Theorem, $X:=\{X_s, t \leq s \leq T\}$ under the new measure $\widehat{\mathbb{P}}$ evolves as  
\begin{align}\label{x}
d X_s =X_s[(\beta+1)(\rho-r+M_s)+\mu_1]ds + X_s \sigma_1 d\widehat{B}_s, \quad X_t=x(=z^{\frac{1}{1-\gamma}}y^{\frac{\gamma}{1-\gamma}}),
\end{align}
where
\refstepcounter{equation}
\begin{align*}
\beta&:=\frac{\gamma}{1-\gamma},\qquad \sigma_1:=\beta\sigma_y-(\beta+1)\theta,\\
\mu_1&:=\frac{1}{2}(\beta+1)\beta\theta^2+\beta\mu_y+\frac{1}{2}(\beta-1)\beta\sigma_y^2-(\beta+1)\beta\theta\sigma_y+\sigma_1(\sigma_y-\theta).\tag{\theequation}\label{Xcoefficients}
\end{align*}
The process \blue{$\widehat B_s:=B_s-B_t-(\sigma_y-\theta)(s-t)$, $s\geq t$, is a Brownian motion under $\widehat{\mathbb P}$.} Then, the solution to (\ref{x}) may be expressed as
\begin{align}\label{xdetailed}
X_s = x \exp \bigg( \int^s_t \Big[(\beta+1)(\rho-r+M_u)+\mu_1-\frac{\sigma_1^2}{2}\Big]du + \int^s_t\sigma_1d \widehat{B}_u  \bigg), \quad \text{for} \ s\geq t,
\end{align}
so that $X$ depends on both initial values $x$ and $m$.

For future frequent use, we also introduce here the new probability measure $\mathbb{\widetilde{P}}$ on $(\Omega, \mathcal{F}_T)$ through the density process
\begin{align}\label{H}
H_s:=\frac{d\mathbb{\widetilde{P}}}{d \widehat{ \mathbb{P}}}\bigg|_{\mathcal{F}_s}= \exp\bigg\{{ \int^s_t \frac{\gamma-1}{\gamma} \sigma_1 d\widehat{ B}_u - \int^s_t \frac{1}{2}\Big(\frac{\gamma-1}{\gamma}\Big)^2 \sigma_1^2 du}\bigg\}, \quad s\in[t,T],
\end{align}
and notice that 
\begin{align}\label{XN}
(X_s)^{\frac{\gamma-1}{\gamma}} = x^{\frac{\gamma-1}{\gamma}}H_s N(s,m),
\end{align}
where 
\begin{align}\label{N}
N(s,m):=  \exp \bigg( \int^s_t {\frac{\gamma-1}{\gamma}} \Big[(\beta+1)(\rho-r+M_u)+\mu_1-\frac{\sigma_1^2}{2}\Big]du + \int^s_t\frac{1}{2}(\frac{\gamma-1}{\gamma})^2 \sigma_1^2 du \bigg).
\end{align}
{Notice that $(H_s)_{s\in[t,T]}$ is a uniformly integrable martingale. By Girsanov's theorem, {$\widetilde B_s:=\widehat B_s-\frac{\gamma-1}{\gamma}\sigma_1(s-t)$, $s\geq t$, is a Brownian motion under $\mathbb{\widetilde P}$.} {Moreover, given that
\begin{align*}
N(s,m)&=\exp\bigg[\bigg(-\frac{\rho-r}{\gamma}+\frac{\gamma-1}{\gamma}\bigg(\mu_1-\frac{\sigma_1^2}{2}\bigg)+\frac{1}{2}\bigg(\frac{\gamma-1}{\gamma}\bigg)^2\sigma_1^2\bigg)(s-t)-\frac{m}{a\gamma}(e^{a(s-t)}-1)\bigg],
\end{align*}
for each $m>0$ there exist $0<\underline L(m)\leq\overline L(m)<\infty$ such that $\underline L(m)\leq N(s,m)\leq\overline L(m)$ for $s\in[t,T]$.}

\subsection{Preliminary properties of the value function}
To study the optimal stopping problem (\ref{4-2}), we find it convenient to introduce the \blue{\textit{Lagrange formulation}}
\begin{align}\label{4-3}
\widehat{J}(t,x,m):=\widetilde{J}(t,x,m)-\widehat{W}(t,x,m)
\end{align}
with
\begin{align}\label{4-4}
\widehat{W}(t,x,m):=Q(x,m)-q_t.
\end{align}

\blue{For an arbitrary $\tau\in\mathcal S$, applying It\^o's formula to $\{e^{-\kappa(s-t)}[Q(X_s,M_s)-q_s],s\in[t,\tau]\}$ and taking conditional expectations give}
\begin{align*}
&\widehat{\mathbb{E}}_{t,x,m}\Big[e^{-\kappa(\tau-t)}\big(Q( X_{\tau}, M_{\tau})-q_\tau \big)\Big]= Q(x,m)-q_t\\&+
\widehat{ \mathbb{E}}_{t,x,m}\bigg[ \int_{t}^{\tau} e^{-\kappa(s-t)}\mathcal{L}\Big(Q( X_s,M_s)-q_s\Big) ds               \bigg ],
\end{align*}
where, for any $F \in C^{1,2,1}(\mathcal{U})$, the second order differential operator $\mathcal{L}$ is such that 
\begin{align*}
{\mathcal{L}}F:=F_t +\frac{1}{2}\sigma_1^2 x^2F_{xx}+[(\beta+1)(\rho-r+m)+\mu_1]xF_x  +a m F_m-\kappa F.
\end{align*}
Combining (\ref{4-2}), (\ref{4-3}) and (\ref{4-4}), we have
 \begin{align}\label{4-13-1}
\widehat{J}(t,x,m)&=\sup_{\tau\in\mathcal S_{t,T}} \widehat{\mathbb{E}}_{t,x,m}\bigg[\int_t^{\tau}e^{-\kappa(s-t)}\Big({U}^\ast_1(X_s)+\mathcal L(Q(X_s,M_s)-q_s)\Big)ds\bigg] \nonumber \\&
=\sup_{t\leq \tau \leq T}  \widehat{\mathbb{E}}_{t,x,m}\bigg[    \int_t^{\tau }  e^{-\kappa(s-t)}\Big( {U}^\ast_1(X_s) -U^\ast_2(X_s)+1 \Big)ds \bigg] \nonumber \\
&=\sup_{t\leq \tau \leq T}  \widehat{\mathbb{E}}_{t,x,m}\bigg[    \int_t^{\tau }  e^{-\kappa(s-t)}\Big( (1-K^{\frac{1-\gamma}{\gamma}}){U}^\ast_1(X_s) +1 \Big)ds \bigg], 
\end{align}
where we have used the fact that 
\begin{align*}
\mathcal{L}(Q(X_s,M_s)-q_s)&= \mathcal{L}Q(X_s,M_s) -\mathcal{L}(q_s)=-U^\ast_2(X_s)+1.
\end{align*}
\blue{Indeed, the homogeneity relation (\ref{4-4-1}) gives
\refstepcounter{equation}
\begin{align*}
Q(z,m)=zyQ(x,m),\qquad x=z^{\frac{1}{1-\gamma}}y^{\frac{\gamma}{1-\gamma}}.\tag{\theequation}\label{wxm}
\end{align*}
}
\blue{Direct substitution gives $\mathcal LQ(x,m)=-U_2^\ast(x)$. Moreover, $q'_s=-e^{-\kappa(T-s)}$, and therefore $\mathcal L(q_s)=q'_s-\kappa q_s=-1$.}

\blue{The augmented Markov family $(X,M)$ is time-homogeneous. More precisely, conditionally on $(X_t,M_t)=(x,m)$, the law of $((X_{t+s},M_{t+s}))_{0\leq s\leq T-t}$ equals the law of $((X_s,M_s))_{0\leq s\leq T-t}$ started from $(X_0,M_0)=(x,m)$. We denote expectation for the latter process by $\widehat{\mathbb E}_{x,m}$. Hence, from (\ref{4-13-1}),}
 \begin{align}\label{widehat{J}}
\widehat{J}(t,x,m)=\sup_{0 \leq \tau \leq T-t} \widehat{\mathbb{E}}_{x,m}\bigg[    \int_0^{\tau }  e^{-\kappa s}\Big( (1-K^{\frac{1-\gamma}{\gamma}}){U}^\ast_1(X_s) +1 \Big)ds \bigg]. 
\end{align}
\blue{Here, the state dynamics are}
\begin{align*}
\blue{M_s=me^{as},\qquad dX_s=X_s[(\beta+1)(\rho-r+M_s)+\mu_1]ds+\sigma_1X_sd\widehat B_s,\qquad X_0=x.}
\end{align*}

As usual in optimal stopping theory, we let
\begin{align}\label{region}
\mathcal{W}:=\{ (t,x,m) \in \mathcal{U} : \widehat{J}(t, x,m)>0  \}, \quad
\mathcal{I}:=\{ (t,x,m) \in  \mathcal{U} : \widehat{J}(t, x,m)=0  \}
\end{align}
be the so-called continuation (waiting) and stopping (retiring) regions, respectively. We denote by $\partial \mathcal{W}$ the boundary of the set $\mathcal{W}.$

Since, for any stopping time $\tau$, the mapping $(x,m) \mapsto \widehat{\mathbb{E}}_{x,m}[    \int_0^{\tau }  e^{-\kappa s}[ (1-K^{\frac{1-\gamma}{\gamma}}){U}^\ast_1(X_s) +1 ]ds ]$ is continuous, $\widehat{J}$ is lower-semicontinuous on $\mathcal{U}$. Hence, $\mathcal{W}$ is open, $\mathcal{I}$ is closed, and the stopping time 
\begin{align}\label{optimalstopping}
\tau^\ast(t,x,m) := \inf\{s \geq 0 :(t+s,X_s,M_s) \in \mathcal{I}\} \wedge (T-t), \quad \widehat{\mathbb{P}}_{x,m}\text{-}{a.s.},
\end{align}
{with the convention $\inf\emptyset=+\infty$, is optimal for $\widehat J(t,x,m)$; see Corollary I.2.9 in \citet{peskir2006optimal}.}

\begin{pro}\label{monotonic}
The value function $\widehat{J}$ is such that $0\leq \widehat{J}(t,x,m)\leq \frac{1}{\kappa}(1-e^{-\kappa(T-t)})$ for all $(t,x,m) \in \mathcal{U}$ and it satisfies the following properties:
\begin{enumerate}[(i)]
	\item  When $0<\gamma<1$, $x \mapsto \widehat{J}(t,x,m) $ is non-decreasing for all $(t,m)  \in [0,T] \times \mathbb{R}_+$;  when $\gamma>1$, $x \mapsto \widehat{J}(t,x,m) $ is non-increasing for all $(t,m)  \in [0,T] \times \mathbb{R}_+$;
	\item $t \mapsto \widehat{J}(t,x,m) $ is non-increasing for all $(x,m)  \in \mathbb{R}^2_+$;
	\item  $m \mapsto \widehat{J}(t,x,m) $ is non-decreasing for all $(t,x)  \in [0,T] \times \mathbb{R}_+.$
\end{enumerate}
\end{pro}

\begin{proof}
The proof is given in Appendix \ref{promonotonic}.
\end{proof}
\blue{For the regularity results below, we impose the following standing nondegeneracy condition.}
\begin{assume}\label{nondegenerate}
\blue{Recall (\ref{Xcoefficients}). We assume $\sigma_1\neq0$, i.e., $\theta\neq\gamma\sigma_y$}.
\end{assume}

The next technical result states some properties of $\widehat{J}$ that will be useful in the study of the regularity of the boundary $\partial \mathcal{W}$. {Here and in the following, $\widetilde{\mathbb E}$ denotes expectation under the probability measure $\widetilde{\mathbb P}$ defined in (\ref{H}).} 
\begin{pro}\label{lip}
Recall $\kappa=r-\mu_y+\sigma_y\theta$, (\ref{H}) and (\ref{N}). {The function $\widehat J$ is locally Lipschitz continuous on $\mathcal U$. At a.e. $(t,x,m)\in\mathcal U$, its weak derivatives satisfy}
\begin{align}\label{lipx}
\widehat{J}_x(t,x,m) &= (K^{\frac{1-\gamma}{\gamma}}-1) x^{-\frac{1}{\gamma}} \widehat{\mathbb{E}}_{x,m}\bigg[    \int_0^{\tau^\ast }  e^{-\kappa s}H_s N(s,m)ds \bigg],
\end{align}
\begin{align}\label{lipm}
\widehat{J}_m(t,x,m)&=  \blue{(K^{\frac{1-\gamma}{\gamma}}-1) \frac{x^{\frac{\gamma-1}{\gamma}}}{1-\gamma}\widehat{\mathbb{E}}_{x,m}\bigg[\int_0^{\tau^\ast}e^{-\kappa s}\frac{e^{as}-1}{a}H_sN(s,m)ds\bigg].}
\end{align}
{Moreover, there exists a constant $C>0$, independent of $(t,x,m)$, such that, at a.e. $(t,x,m)\in\mathcal U$,}
\begin{align}\label{lipt}
-\frac{1}{T-t}( C x^{\frac{\gamma-1}{\gamma}}	\widetilde{\mathbb{E}}[\tau^\ast]+	\widehat{\mathbb{E}}[\tau^\ast]) \leq \widehat{J}_t(t,x,m) \leq 0,
\end{align}
where $\tau^\ast:=\tau^\ast(t, x, m)$ is the optimal stopping time for the problem with initial data $(t, x, m)$. 
\end{pro}
\begin{proof}
The proof is given in Appendix \ref{prolip}.

%\textcolor{red}{what is the main difference from Theorem 3.1 of \citet{de2019lipschitz}? }

\end{proof}

We conclude with asymptotic limits of $\widehat{J}$.
\begin{pro}\label{limit}
	When $\gamma<1$ we have
	\begin{align*}
		\lim_{x \to 0}\widehat{J}(t,x,m) =0, \quad \lim_{x\to \infty}\widehat{J}(t,x,m)=\frac{1}{\kappa}\Big(1-e^{-\kappa(T-t)}\Big), \ \text{for all} \ (t,m)\in [0,T) \times \mathbb{R}_+;
		\end{align*}  
		 when $\gamma>1$ we have 
		\begin{align*}
			\lim_{x \to 0}\widehat{J}(t,x,m)=\frac{1}{\kappa}\Big(1-e^{-\kappa(T-t)}\Big), \quad \lim_{x \to \infty}\widehat{J}(t,x,m)=0, \ \text{for all} \ (t,m)\in [0,T) \times \mathbb{R}_+.
					\end{align*}
	\end{pro}

\begin{proof}
The proof is given in Appendix \ref{prolimit}.

\end{proof}

\subsection{Preliminary properties of the optimal boundary}

In this section, we show that the boundary $\partial \mathcal{W}$ can be represented by a function $b(t,m)$. We establish connectedness of the sets $\mathcal{W}$ and $\mathcal{I}$ with respect to the $x$-variable and give some preliminary properties of the optimal boundary. 

\begin{lemma}\label{b}
 \blue{Let $L:=\big[(K^{\frac{1-\gamma}{\gamma}}-1)\frac{\gamma}{1-\gamma}\big]^{\frac{\gamma}{1-\gamma}}>0$.} \blue{When $\gamma<1$, there exists a boundary $b:[0,T)\times\mathbb R_+\to[0,L]$ such that}
\begin{align}\label{stopping1}
\blue{\mathcal{I} \cap \{t<T\} =\{(t,x,m)\in[0,T)\times\mathbb R_+^2:0<x\leq b(t,m)\}.}
\end{align}
Moreover, the function $b$ has the following properties:
\begin{enumerate}[(i)]
	\item $t \to b(t,m)$ is non-decreasing for any $m \in \mathbb{R}_+;$
\item $m \to b(t,m)$ is non-increasing for any $t \in [0,T)$;
	\item \blue{one has $0\leq b(t,m)\leq L$.}
\end{enumerate}

\blue{When $\gamma>1$, there exists a boundary $b:[0,T)\times\mathbb R_+\to[L,\infty]$ such that}
\begin{align}\label{stopping2}
\blue{\mathcal{I} \cap \{t<T\} =\{(t,x,m)\in[0,T)\times\mathbb R_+^2:x\geq b(t,m)\}.}
\end{align}

Moreover, the function $b$ has the following properties:
\begin{enumerate}[(i)]
	\item $t \to b(t,m)$ is non-increasing for any $m \in \mathbb{R}_+$;
\item $m \to b(t,m)$ is non-decreasing for any $t \in [0,T)$;
	\item  \blue{one has $b(t,m)\geq L$.}
\end{enumerate}
\end{lemma}

\begin{proof}
The proof is given in Appendix \ref{prob}.
\end{proof}

\subsection{Characterization of the free boundary and the value function}
\blue{To establish the value-function regularity needed to obtain an integral equation for the optimal boundary, we now perform a change of variable.}
{\color{black}
{\color{black}
Fix $m>0$ and set $\overline m(t):=me^{at}$. The fixed initial level $m$ is suppressed in the following notation. Define
\refstepcounter{equation}
\begin{align*}
\mathcal J(t,x)&:=\widehat J(t,x,\overline m(t))\\
&=\sup_{0\leq\tau\leq T-t}\widehat{\mathbb E}_{t,x}\bigg[\int_0^\tau e^{-\kappa s}\Big((1-K^{\frac{1-\gamma}{\gamma}})U_1^\ast(X_s)+1\Big)ds\bigg],\tag{\theequation}\label{mathcalJ}
\end{align*}
where
\begin{align*}
dX_s&=X_s\Big[(\beta+1)(\rho-r+\overline m(t+s))+\mu_1\Big]ds+\sigma_1X_sd\widehat B_s,\qquad X_0=x.
\end{align*}
}
For $F\in C^{1,2}([0,T]\times\mathbb R_+)$, let
\begin{align*}
\overline{\mathcal L}F:=F_t+\frac{1}{2}\sigma_1^2x^2F_{xx}+\Big[(\beta+1)(\rho-r+\overline m(t))+\mu_1\Big]xF_x-\kappa F.
\end{align*}
Set $\overline{\mathcal W}:=\{(t,x)\in[0,T)\times\mathbb R_+:\mathcal J(t,x)>0\},$ and 
$\overline{\mathcal I}:=\{(t,x)\in[0,T]\times\mathbb R_+:\mathcal J(t,x)=0\}.$ The optimal stopping time is 
\begin{align*}
\overline\tau^\ast(t,x)&:=\inf\{s\in[0,T-t]:(t+s,X_s)\in\overline{\mathcal I}\}.
\end{align*}
Clearly, $\tau^\ast(t,x,\overline{m}(t))=\overline{\tau}^\ast(t,x)$ a.s., and the corresponding boundary is 
\begin{align*}
	\overline b(t)&:=b(t,\overline m(t)),\qquad t\in[0,T).
\end{align*}

\begin{theorem}\label{lipschitz}
For every $m>0$, the boundary $\overline b$ is finite, strictly positive, and locally Lipschitz continuous on $[0,T)$. In particular, $0<b(t,m)<\infty$ for every $(t,m)\in[0,T)\times\mathbb R_+$.
\end{theorem}

\begin{proof}
The proof is given in Appendix \ref{prolipschitz}.
\end{proof}
}

{\color{black}
\begin{lemma}\label{interior}
When $\gamma<1$, set
\begin{align*}
\widehat\tau(t,x)&:=\inf\{s\geq0:X_s<\overline b(t+s)\}\wedge(T-t),\\
\overline\tau^\ast(t,x)&:=\inf\{s\geq0:X_s\leq\overline b(t+s)\}\wedge(T-t).
\end{align*}
Then $\widehat\tau(t,x)=\overline\tau^\ast(t,x)$ a.s. When $\gamma>1$, the same conclusion holds with
\begin{align*}
\widehat\tau(t,x)&:=\inf\{s\geq0:X_s>\overline b(t+s)\}\wedge(T-t),\\
\overline\tau^\ast(t,x)&:=\inf\{s\geq0:X_s\geq\overline b(t+s)\}\wedge(T-t).
\end{align*}
\end{lemma}

\begin{proof}
The proof is given in Appendix \ref{prointerior}.
\end{proof}

The previous lemma in turn yields the following continuity property of $\overline{\tau}^\ast$, which will then be fundamental in the proof of Proposition \ref{regularity} below. 

\begin{pro}\label{stopping}
For every sequence $(t_n,x_n)\to(t,x)$ in $[0,T]\times\mathbb R_+$, one has
\begin{align*}
\overline\tau^\ast(t_n,x_n)\longrightarrow\overline\tau^\ast(t,x),\qquad\widehat{\mathbb P}\text{-a.s.}
\end{align*}
\end{pro}

\begin{proof}
The proof exploits Lemma \ref{interior} and arguments completely analogous to those employed in the proof of Proposition 5.2 in \citet{de2017dividend} and it is therefore omitted.
\end{proof}

\begin{pro}\label{regularity}
The value function satisfies $\mathcal J\in C^{1,2}(\overline{\mathcal W})\cap C^1([0,T)\times\mathbb R_+)\cap C([0,T]\times\mathbb R_+),$ and $
\mathcal J_{xx}\in L^\infty_{loc}([0,T)\times\mathbb R_+).
$
Moreover,
{\color{black}
\refstepcounter{equation}
\begin{align*}
\begin{cases}
(\overline{\mathcal L}\mathcal J)(t,x)=-(1-K^{\frac{1-\gamma}{\gamma}})U_1^\ast(x)-1, &(t,x)\in\overline{\mathcal W},\\
\mathcal J(t,x)=0, &(t,x)\in\overline{\mathcal I},\\
\mathcal J(T,x)=0, &x\in\mathbb R_+,\\
\mathcal J_t(t,x)=\mathcal J_x(t,x)=0, &(t,x)\in\partial\overline{\mathcal W}\cap\{t<T\}.
\end{cases}\tag{\theequation}\label{classcial}
\end{align*}
}
\end{pro}

\begin{proof}
On every parabolic cylinder compactly contained in $\overline{\mathcal W}$, {$x$ is bounded away from zero.} Assumption \ref{nondegenerate} makes $\overline{\mathcal L}$ uniformly parabolic there. Standard interior parabolic regularity, as in Corollary 2.4.3 of \citet{krylov2008lectures} and Theorem 2.7.7 of \citet{karatzas1998methods}, gives $\mathcal J\in C^{1,2}(\overline{\mathcal W})$ and the first equation in (\ref{classcial}).

{\color{black}
It remains to study regularity across the boundary. Proposition \ref{lip} shows that $\widehat J$ is locally Lipschitz, with weak derivatives given by (\ref{lipx})--(\ref{lipm}) and the bound (\ref{lipt}). Hence $\mathcal J$ is locally Lipschitz as well. Since $\mathcal J\in C^{1,2}(\overline{\mathcal W})$, its weak derivatives are classical derivatives in $\overline{\mathcal W}$. 

  Now fix $(t_0,x_0)\in\partial\overline{\mathcal W}$ with $t_0<T$ and take $(t_n,x_n)\in\overline{\mathcal W}$ converging to $(t_0,x_0)$. Proposition \ref{stopping} gives $\overline\tau^\ast(t_n,x_n)\to0$ a.s. Using (\ref{lipm}), the locally uniform bounds for $N$ and $(e^{as}-1)/a$, and the density $H$, we obtain
\begin{align*}
0\leq\widehat J_m(t_n,x_n,\overline m(t_n))
\leq C\widetilde{\mathbb E}[\overline\tau^\ast(t_n,x_n)]\longrightarrow0.
\end{align*}
Here $C$ is locally uniform. (\ref{lipx}) similarly gives $\mathcal J_x(t_n,x_n)\to0$. From (\ref{lipt}), with $\tau:=\overline\tau^\ast(t,x)$,
\begin{align*}
-\frac{Cx^{\frac{\gamma-1}{\gamma}}\widetilde{\mathbb E}[\tau]+\widehat{\mathbb E}[\tau]}{T-t}
&\leq\mathcal J_t(t,x)-a\overline m(t)\widehat J_m(t,x,\overline m(t))\leq0.
\end{align*}
The time estimate, bounded convergence for the stopping-time expectations, and $t_0<T$ then yield $\mathcal J_t(t_n,x_n)\to0$. 
 Both derivatives vanish in the interior of $\overline{\mathcal I}$, so local Lipschitz continuity and these continuous boundary limits give the asserted $C^1$ extension and smooth fit.
}

Solving the first equation in (\ref{classcial}) for the second derivative gives, in $\overline{\mathcal W}$,
\begin{align*}
\frac{1}{2}\sigma_1^2x^2\mathcal J_{xx}=-(1-K^{\frac{1-\gamma}{\gamma}})U_1^\ast(x)-1-\mathcal J_t
-\Big[(\beta+1)(\rho-r+\overline m(t))+\mu_1\Big]x\mathcal J_x+\kappa\mathcal J.
\end{align*}
The right-hand side is bounded on every compact subset of $[0,T)\times\mathbb R_+$. Define $\mathcal J_{xx}=0$ in the interior of $\overline{\mathcal I}$. Since $\overline b$ is locally Lipschitz and $\mathcal J_x$ is continuous across its graph, the distributional derivative has no singular boundary term. Hence the piecewise-defined function is the weak second derivative and belongs to $L^\infty_{loc}$. Finally,
\begin{align*}
0\leq\mathcal J(t,x)\leq\int_0^{T-t}e^{-\kappa s}ds\longrightarrow0
\end{align*}
as $t\uparrow T$, which proves continuity at the terminal time.
\end{proof}

We are now in the conditions of determining a nonlinear integral equation that characterizes uniquely the free boundary. Such a characterization results from an integral representation of the value function $\mathcal{J}$. This is accomplished by the next theorem, which exploits the regularity properties of $\mathcal{J}$ proved so far.

{\color{black}
\begin{theorem}\label{integral}
If $\gamma<1$, the value function in (\ref{mathcalJ}) has the following representation for every $(t,x)\in[0,T)\times\mathbb R_+$:
\refstepcounter{equation}
\begin{align*}
\mathcal J(t,x)=\widehat{\mathbb E}_{t,x}\bigg[\int_0^{T-t}e^{-\kappa s}\Big((1-K^{\frac{1-\gamma}{\gamma}})U_1^\ast(X_s)+1\Big)\mathds{1}_{\{X_s\geq\overline b(t+s)\}}ds\bigg].\tag{\theequation}\label{integralJ1}
\end{align*}
Moreover, with $L$ as in Lemma \ref{b}, $\overline b$ is the unique continuous function on $[0,T]$, strictly positive and not exceeding $L$, that satisfies
\refstepcounter{equation}
\begin{align*}
0=\widehat{\mathbb E}_{t,\overline b(t)}\bigg[\int_0^{T-t}e^{-\kappa s}\Big((1-K^{\frac{1-\gamma}{\gamma}})U_1^\ast(X_s)+1\Big)\mathds{1}_{\{X_s\geq\overline b(t+s)\}}ds\bigg],\tag{\theequation}\label{integralb1}
\end{align*}
for every $t<T$, with terminal condition $\overline b(T)=L$.

If $\gamma>1$, the corresponding representation is
\refstepcounter{equation}
\begin{align*}
\mathcal J(t,x)=\widehat{\mathbb E}_{t,x}\bigg[\int_0^{T-t}e^{-\kappa s}\Big((1-K^{\frac{1-\gamma}{\gamma}})U_1^\ast(X_s)+1\Big)\mathds{1}_{\{X_s\leq\overline b(t+s)\}}ds\bigg].\tag{\theequation}\label{integralJ2}
\end{align*}
In this case, $\overline b$ is the unique continuous, finite function on $[0,T]$, not below $L$, that satisfies
\refstepcounter{equation}
\begin{align*}
0=\widehat{\mathbb E}_{t,\overline b(t)}\bigg[\int_0^{T-t}e^{-\kappa s}\Big((1-K^{\frac{1-\gamma}{\gamma}})U_1^\ast(X_s)+1\Big)\mathds{1}_{\{X_s\leq\overline b(t+s)\}}ds\bigg],\tag{\theequation}\label{integralb2}
\end{align*}
for every $t<T$, with terminal condition $\overline b(T)=L$.
\end{theorem}

\begin{proof}
\textbf{Step 1:} We prove (\ref{integralJ1}); the case $\gamma>1$ is analogous. Fix $(t,x)\in[0,T)\times\mathbb R_+$ and, for all sufficiently large $n$, set
\begin{align*}
\tau_n:=\inf\{s\geq0:X_s\notin(1/n,n)\}\wedge(T-t-1/n).
\end{align*}
Continuity and strict positivity of $X$ give $\tau_n\uparrow T-t$ a.s. By Proposition \ref{regularity}, $\mathcal J\in C^1([0,T)\times\mathbb R_+)$ and its weak second spatial derivative is locally bounded. Moreover, $X_s$ has a density for every $s>0$, so
\begin{align*}
\widehat{\mathbb E}_{t,x}\bigg[\int_0^{T-t}\mathds{1}_{\{X_s=\overline b(t+s)\}}ds\bigg]=0.
\end{align*}
 The weak Dynkin formula in Lemma 8.1 and Theorem 8.5 of \citet{bensoussan2011applications}, pp.\ 183--186, therefore gives
\begin{align*}
\mathcal J(t,x)=\widehat{\mathbb E}_{t,x}\bigg[e^{-\kappa\tau_n}\mathcal J(t+\tau_n,X_{\tau_n})-\int_0^{\tau_n}e^{-\kappa s}(\overline{\mathcal L}\mathcal J)(t+s,X_s)ds\bigg].
\end{align*}
Using (\ref{classcial}) in the continuation region and $\mathcal J=0$ in the interior of the stopping region, we obtain
\begin{align*}
\mathcal J(t,x)=\widehat{\mathbb E}_{t,x}\bigg[e^{-\kappa\tau_n}\mathcal J(t+\tau_n,X_{\tau_n})
+\int_0^{\tau_n}e^{-\kappa s}\Big((1-K^{\frac{1-\gamma}{\gamma}})U_1^\ast(X_s)+1\Big)\mathds{1}_{\{X_s\geq\overline b(t+s)\}}ds\bigg].
\end{align*}
Notice that $(1-K^{\frac{1-\gamma}{\gamma}})U_1^\ast(X_s)+1$ is bounded by $1+CX_s^{\frac{\gamma-1}{\gamma}}$, whose time integral is integrable by the lognormal moment bounds. Also, Proposition \ref{monotonic} gives
\begin{align*}
0\leq\mathcal J(t+\tau_n,X_{\tau_n})\leq T-t-\tau_n\longrightarrow0.
\end{align*}
Dominated convergence now proves (\ref{integralJ1}). Reversing the indicator inequality gives (\ref{integralJ2}).

\textbf{Step 2:} We follow the argument in Proposition 4.10 of \citet{de2019free}. Here, however, monotonicity of $\overline b$ is not available, so we first establish existence of its terminal limit by comparing current mortality levels. In this step, set $\mathcal C:=\frac{e^{aT}-1}{a|1-\gamma|}$. For $m>0$ and $\delta\geq0$, (\ref{Nm}) gives
\begin{align*}
N(s,m+\delta)=N(s,m)\exp\bigg(-\frac{\delta}{\gamma}\frac{e^{as}-1}{a}\bigg).
\end{align*}
When $\gamma<1$, for $s\leq T-t$,
\begin{align*}
\big(xe^{-\mathcal C\delta}\big)^{\frac{\gamma-1}{\gamma}}N(s,m+\delta)
&=x^{\frac{\gamma-1}{\gamma}}N(s,m)\exp\bigg(\frac{\delta}{\gamma}\frac{e^{aT}-e^{as}}{a}\bigg)
\geq x^{\frac{\gamma-1}{\gamma}}N(s,m).
\end{align*}
The same inequality holds with $xe^{\mathcal C\delta}$ when $\gamma>1$. Since the coefficient of $x^{\frac{\gamma-1}{\gamma}}H_sN(s,m)$ in the running payoff is negative, comparison under every common stopping time and Proposition \ref{monotonic} give
\begin{align*}
\widehat J(t,xe^{-\mathcal C\delta},m+\delta)&\leq\widehat J(t,x,m)\leq\widehat J(t,x,m+\delta),\qquad\gamma<1,\\
\widehat J(t,xe^{\mathcal C\delta},m+\delta)&\leq\widehat J(t,x,m)\leq\widehat J(t,x,m+\delta),\qquad\gamma>1.
\end{align*}
Moreover, from Lemma \ref{b}, we obtain
\begin{align*}
e^{-\mathcal C\delta}b(t,m)&\leq b(t,m+\delta)\leq b(t,m),\qquad\gamma<1,\\
b(t,m)&\leq b(t,m+\delta)\leq e^{\mathcal C\delta}b(t,m),\qquad\gamma>1.
\end{align*}
Fix the mortality level $\overline m(T)$. By Lemma \ref{b}, the function $t\mapsto b(t,\overline m(T))$ is nondecreasing and bounded above by $L$ when $\gamma<1$, and nonincreasing and bounded below by $L$ when $\gamma>1$. Since Theorem \ref{lipschitz} gives $0<b(0,\overline m(T))<\infty$, this function has a finite, strictly positive limit as $t\uparrow T$. Taking $m=\overline m(t)$ and $\delta=\overline m(T)-\overline m(t)$ in the preceding comparison gives
\begin{align*}
b(t,\overline m(T))&\leq\overline b(t)\leq e^{\mathcal C(\overline m(T)-\overline m(t))}b(t,\overline m(T)),\qquad\gamma<1,\\
e^{-\mathcal C(\overline m(T)-\overline m(t))}b(t,\overline m(T))&\leq\overline b(t)\leq b(t,\overline m(T)),\qquad\gamma>1.
\end{align*}
The exponential factors converge to one. The squeeze theorem therefore shows that $\overline b(T-)$ exists and equals the finite, strictly positive limit of $b(t,\overline m(T))$.

Lemma \ref{b} gives $\overline b(T-)\leq L$ when $\gamma<1$ and $\overline b(T-)\geq L$ when $\gamma>1$. Suppose that the relevant inequality is strict. Choose $\overline b(T-)<x_1<x_2<L$ when $\gamma<1$, and $L<x_1<x_2<\overline b(T-)$ when $\gamma>1$. In either case $(s,x)\in\overline{\mathcal W}$ for $x\in(x_1,x_2)$ and all $s<T$ sufficiently close to $T$. Let $\varphi\in C_c^\infty(x_1,x_2)$ satisfy $\varphi\geq0$ and $\int_{x_1}^{x_2}\varphi(x)dx=1$. The PDE in (\ref{classcial}) and integration by parts give
\begin{align*}
\int_{x_1}^{x_2}\mathcal J_t(s,x)\varphi(x)dx
&=\int_{x_1}^{x_2}\Big((K^{\frac{1-\gamma}{\gamma}}-1)U_1^\ast(x)-1\Big)\varphi(x)dx\\
&\quad+\int_{x_1}^{x_2}\mathcal J(s,x)\bigg[\kappa\varphi(x)-\frac{\sigma_1^2}{2}(x^2\varphi(x))''\\
&\hspace{85pt}+\Big((\beta+1)(\rho-r+\overline m(s))+\mu_1\Big)(x\varphi(x))'\bigg]dx.
\end{align*}
The first integral is strictly positive by the choice of $(x_1,x_2)$. The second tends to zero as $s\uparrow T$, because $0\leq\mathcal J(s,x)\leq T-s$ and the coefficients are bounded on this interval. Thus $\int_{x_1}^{x_2}\mathcal J_t(s,x)\varphi(x)dx>0$ for $s$ sufficiently close to $T$. Integrating up to $T$ gives, for sufficiently small $\delta>0$,
\begin{align*}
0<\int_{T-\delta}^T\int_{x_1}^{x_2}\mathcal J_t(s,x)\varphi(x)dxds
=-\int_{x_1}^{x_2}\mathcal J(T-\delta,x)\varphi(x)dx<0,
\end{align*}
a contradiction.  Therefore $\overline b(T-)=L$, and we set $\overline b(T):=L$.

\textbf{Step 3:} Substituting $x=\overline b(t)$ into (\ref{integralJ1}) or (\ref{integralJ2}) and using value matching gives (\ref{integralb1}) or (\ref{integralb2}). Uniqueness in the stated class follows from the four-step strong-Markov argument in Theorem 3.1 of \citet{peskir2005american}. Since the present setting does not create additional difficulties we omit further details.
\end{proof}

\begin{remark}
The representations in $(t,x,m)$ and $(t,x)$ serve different purposes. When working with $\widehat{J}(t,x,m)$, $(X,M)$ is time-homogeneous, so varying $t$ at fixed $(x,m)$ only changes the remaining horizon. This gives the time monotonicity of $\widehat J$ in Proposition \ref{monotonic}. This also yields the mortality and time monotonicity of the boundary $b(t,m)$ in Lemma \ref{b}. Once a mortality path is fixed, its deterministic component is absorbed into time. The resulting $\mathcal J(t,x)$ problem permits one-dimensional interior regularity and standard parabolic estimates.

However, the monotonicity of $b(t,m)$ does not imply monotonicity of $\overline b(t)=b(t,\overline m(t))$. Theorem \ref{integral} does not assume that $\overline b$ is monotone. Its terminal limit follows from the mortality comparison with $b(t,\overline m(T))$ and the squeeze theorem in Step 2. Thus the $\widehat{J}(t,x,m)$ formulation supplies useful comparisons, while the $\mathcal J(t,x)$ formulation supplies the regularity needed for the integral characterization and the subsequent duality argument.
\end{remark}
}
}

\section{Optimal boundaries and strategies in terms of primal variables}\label{sec5}
In the previous section, we studied the properties of the dual value function $J(t,z,m,y)$ in the coordinates $(t,z,m,y)$, where $t$ denotes time, $z$ denotes marginal utility, $m$ denotes the force of mortality, and $y$ denotes labor income. We now return to the primal value function $V(t,w,m,y)$, where $w$ denotes the agent's financial wealth.

\begin{lemma}\label{convex}
For every fixed $(t,m,y)$, the function $z\mapsto J(t,z,m,y)$ belongs to $C^1(\mathbb R_+)$ and is strictly decreasing and strictly convex on $\mathbb R_+$. Moreover, 
\begin{align*}
	\lim_{z\to 0_+}-J_z(t,z,m,y)=+\infty \quad \text{and} \quad \lim_{z \to +\infty} -J_z(t,z,m,y)=0.
\end{align*}
\end{lemma}
\begin{proof}

Fix $(t,m,y)$. At $t=T$, the result follows directly from $J(T,z,m,y)=Q(z,m)$, so assume $t<T$. Set $x=z^{\frac{1}{1-\gamma}}y^{\frac{\gamma}{1-\gamma}}$. By (\ref{jandj}), (\ref{4-3}), (\ref{4-4}), and (\ref{wxm}),
\begin{align}\label{r1}
J(t,z,m,y)=Q(z,m)-zg_t+zy\widehat J(t,x,m).
\end{align}
\blue{Along the mortality path passing through $(t,m)$, one has $\widehat J(t,x,m)=\mathcal J(t,x)$. Proposition \ref{regularity} gives $\mathcal J\in C^1$, so $x\mapsto\widehat J(t,x,m)$ is $C^1$ and $\widehat J_x(t,x,m)=\mathcal J_x(t,x)$. Equation (\ref{r1}) and the smooth dependence of $x$ on $z$ therefore imply $J(t,\cdot,m,y)\in C^1(\mathbb R_+)$. In particular,}
\begin{align}\label{Jzformula}
J_z(t,z,m,y)=Q_z(z,m)-g_t+y\widehat J(t,x,m)+\frac{xy}{1-\gamma}\widehat J_x(t,x,m).
\end{align}

Next, for $\tau\in\mathcal S$, define
\begin{align*}
\overline J(t,z,m,y;\tau):=\mathbb E_{t,z,m,y}\bigg[\int_t^\tau e^{-\int_t^s(\rho+M_u)du}U_1^\ast(Z_s)ds+e^{-\int_t^\tau(\rho+M_u)du}\big(Q(Z_\tau,M_\tau)-Z_\tau g_\tau\big)\bigg].
\end{align*}
Then $J(t,z,m,y)=\sup_{\tau\in\mathcal S}\overline J(t,z,m,y;\tau)$.

 We now establish strict monotonicity and strict convexity by adapting the argument of Lemma 8.1 in \citet{karatzas2000utility}. Since $Z_s=zP_s$, the strict monotonicity and strict convexity of $U_1^\ast$ and $Q(\cdot,m)$, together with the linearity of $-Z_\tau g_\tau$, show that $z\mapsto\overline J(t,z,m,y;\tau)$ is strictly decreasing and strictly convex for every fixed $\tau$. If $0<z_1<z_2$ and $\tau_2$ is optimal for $J(t,z_2,m,y)$, then
\begin{align*}
J(t,z_1,m,y)\geq\overline J(t,z_1,m,y;\tau_2)>\overline J(t,z_2,m,y;\tau_2)=J(t,z_2,m,y).
\end{align*}
Likewise, if $z_\alpha:=\alpha z_1+(1-\alpha)z_2$, with $\alpha\in(0,1)$, and $\tau_\alpha$ is optimal for $J(t,z_\alpha,m,y)$, then
\begin{align*}
J(t,z_\alpha,m,y)&=\overline J(t,z_\alpha,m,y;\tau_\alpha)\\
&<\alpha\overline J(t,z_1,m,y;\tau_\alpha)+(1-\alpha)\overline J(t,z_2,m,y;\tau_\alpha)\\
&\leq\alpha J(t,z_1,m,y)+(1-\alpha)J(t,z_2,m,y).
\end{align*}
Hence, $J(t,\cdot,m,y)$ is strictly decreasing and strictly convex.

\color{black}
We next prove the limiting behavior. From Proposition \ref{monotonic},
we have $0\leq\widehat J(t,x,m)\leq q_t$. Since $g_t=yq_t$ and (\ref{r1}), it follows that
\begin{align*}
Q(z,m)-zg_t\leq J(t,z,m,y)\leq Q(z,m).
\end{align*}
By (\ref{w}), $Q(z,m)=C(m)z^\frac{\gamma-1}{\gamma}$, with $C(m)>0$ when $\gamma<1$ and $C(m)<0$ when $\gamma>1$. If $\gamma<1$, then $Q(z,m)\to\infty$ as $z\downarrow0$, and the lower bound gives $J(t,z,m,y)\to\infty$. Convexity therefore yields, for any fixed $z_0>0$,
\begin{align*}
J_z(t,z,m,y)\leq\frac{J(t,z_0,m,y)-J(t,z,m,y)}{z_0-z}\longrightarrow-\infty.
\end{align*}
If $\gamma>1$, then $\frac{\gamma-1}{\gamma} \in(0,1)$ and convexity gives
\begin{align*}
J_z(t,z,m,y)&\leq\frac{J(t,2z,m,y)-J(t,z,m,y)}{z}\\
&\leq C(m)(2^\frac{\gamma-1}{\gamma}-1)z^{\frac{-1}{\gamma}}+g_t\longrightarrow-\infty.
\end{align*}
Thus, $-J_z(t,z,m,y)\to\infty$ as $z\downarrow0$ in both cases.

Finally, as $z\uparrow\infty$, Proposition \ref{limit} gives $\widehat J(t,x,m)\to q_t$: for $\gamma<1$ one has $x\uparrow\infty$, whereas for $\gamma>1$ one has $x\downarrow0$. Since $Q(z,m)/z\to0$, (\ref{r1}) yields
\begin{align}\label{r2}
\lim_{z\uparrow\infty}\frac{J(t,z,m,y)}{z}=0.
\end{align}
Because $J(t,\cdot,m,y)$ is convex and decreasing, $J_z(t,\cdot,m,y)$ is nondecreasing and nonpositive. Hence, 
\begin{align*}
	\ell :=\lim_{z\uparrow \infty}J_z(t,z,m,y)
\end{align*}
\blue{exists in $(-\infty,0]$. If $\ell<0$, monotonicity gives $J_z(t,z,m,y)\leq\ell$ for every $z>0$. Integrating from any fixed $z_0>0$ yields $J(t,z,m,y)\leq J(t,z_0,m,y)+\ell(z-z_0)$ for $z\geq z_0$, and hence}
{\color{black}
\begin{align*}
\limsup_{z\uparrow\infty}\frac{J(t,z,m,y)}{z}\leq\ell<0,
\end{align*}
}
contradicting (\ref{r2}). Hence, $J_z(t,z,m,y)\to0$ as $z\uparrow\infty$.
\end{proof}

We now provide the following theorem that establishes a dual relation between the original problem (\ref{2-6}) and the optimal stopping problem (\ref{3-13}) and the corresponding optimal strategies.

%\begin{lemma}
%	Define the primal map $\mathcal{P}$ as 
%	\begin{align*}
%		\mathcal{P}(t,z,m,y)=-J_z(t,z,m,y),
%	\end{align*}	
%	and define the primal continuation region as 
%	\begin{align*}
%		\mathcal{W}^p:= \{(t,w,m,y)\in \mathcal{O}: V(t,w,m,y)> \widehat{V}(w,m)\}.
%	\end{align*}
%	Then we have $\mathcal{W}=\mathcal{W}^p$.
%\end{lemma}
%\begin{proof}
%	Based on duality relation in Theorem \ref{dualityrelation2} and strict convexity of $J$, we have 
%	\begin{align*}
%		V(t,w,m,y)= J(t,z,m,y)+z(w+g_t)
%	\end{align*} 
%	with $z= \mathcal{P}^{-1}(w+g_t)$.
	
%	If $z \in \mathcal{W}$, then 
%	\begin{align*}
%		J(t,z^\ast,m,y)+z^\ast(w+g_t) >+z^\ast(w+g_t)+ 			\end{align*}
%\end{proof}

\begin{theorem}\label{verification}
For given $(t,w,m,y)\in \mathcal O' $, the following duality relations hold:
\begin{align*}
	V(t,w,m,y)=\inf_{z>0}[J(t,z,m,y)+z(w+g_t)], \quad J(t,z,m,y)= \sup_{w>-g_t}[V(t,w,m,y)-z(w+g_t)]. 
\end{align*}
Moreover, let $z^\ast:=z^\ast(t,w,m,y)$ be the unique solution of $w+g_t+J_z(t,z^\ast,m,y)=0$, let $Z^\ast$ solve (\ref{3-7}) with $Z_t^\ast=z^\ast$, and let $\tau^\ast\in\mathcal S$ be the corresponding optimal stopping time for (\ref{3-13}). Before retirement, for $t\leq s<\tau^\ast$, the optimal controls satisfy
\begin{align*}
c_s^\ast=(Z_s^\ast)^{-\frac{1}{\gamma}},\qquad
\pi_s^\ast=\frac{\theta Z_s^\ast J_{zz}(s,Z_s^\ast,M_s,Y_s)-\sigma_yY_s\big(q_s+J_{zy}(s,Z_s^\ast,M_s,Y_s)\big)}{\sigma},
\end{align*}
and the corresponding optimal wealth process is $W_s^\ast=-J_z(s,Z_s^\ast,M_s,Y_s)-g_s$. 
\blue{After $\tau^\ast$, the optimal controls are the post-retirement controls in Theorem \ref{dualityrelation}, initialized at $(\tau^\ast,W_{\tau^\ast}^\ast,M_{\tau^\ast})$.}

\end{theorem}

\begin{proof}
The proof is given in Appendix \ref{proofsolution}.
\end{proof}

From Theorem \ref{verification}, for any $(t,w,m,y) \in \mathcal{O}'$,
\begin{align*}
V(t,w,m,y)=\inf_{z>0}[J(t,z,m,y)+z(w+g_t)].
\end{align*}
\blue{For each fixed $(t,m,y)$, Lemma \ref{convex} shows that $z\mapsto-J_z(t,z,m,y)$ is continuous, strictly decreasing, and maps $(0,\infty)$ onto $(0,\infty)$. Hence, for every $w>-g_t$, there exists a unique $z^\ast:=z^\ast(t,w,m,y)>0$ satisfying $w+g_t=-J_z(t,z^\ast,m,y)$. The infimum is attained at this point, and
\refstepcounter{equation}
\begin{align*}
V(t,w,m,y)=J(t,z^\ast,m,y)+z^\ast(w+g_t).\tag{\theequation}\label{5-1}
\end{align*}
Moreover, for each fixed $(t,m,y)$, the map $w\mapsto z^\ast(t,w,m,y)$ is continuous and strictly decreasing from $(-g_t,\infty)$ onto $(0,\infty)$. Its inverse is $w^\ast(t,z,m,y):=-J_z(t,z,m,y)-g_t$, which is continuous and strictly decreasing from $(0,\infty)$ onto $(-g_t,\infty)$.}

%Since $z \mapsto {J}(t,z,m,y)+z(w+g_t)$ is strictly convex (cf.\ Lemma \ref{convex}), there exists a unique $z^\ast(t,w,m,y)>0$ such that 
%\begin{align}\label{5-1}
%V(t,w,m,y)=J(t,z^\ast(t,w,m,y),m,y)+ z^\ast(t,w,m,y)(w+g_t),
%\end{align} 
%where $z^\ast(t,w,m,y):=\mathcal{I}^{J}(t,-w-g_t,m,y)$, with $\mathcal{I}^J$ being the inverse function of $J_z$. Moreover, $z^\ast\in C(\mathcal{O}')$ and, for any $(t,m,y)\in\mathcal U$, the map $w\mapsto z^\ast(t,w,m,y)$ is strictly decreasing. The two limits in Lemma \ref{convex} show that this map sends $(-g_t,+\infty)$ onto $\mathbb R_+$. It is therefore a bijection and has an inverse $w^\ast(t,\cdot,m,y)$, which is continuous, strictly decreasing, and maps $\mathbb R_+$ onto $(-g_t,+\infty)$.

Let us now define
\begin{equation}\label{newb}
\left\{
\begin{aligned}
\widehat{b}(t,m,y)&:=w^\ast(t,b(t,m)^{1-\gamma}y^{-\gamma},m,y),\\
\mathcal{W}_z &:=\{(t,z,m,y) \in \mathcal{O}: (t, z^{\frac{1}{1-\gamma}}y^{\frac{\gamma}{1-\gamma}}, m) \in \mathcal{W} \},\\
\mathcal{I}_z &:=\{(t,z,m,y) \in  \mathcal{O}:  (t, z^{\frac{1}{1-\gamma}}y^{\frac{\gamma}{1-\gamma}}, m) \in \mathcal{I} \},\\
\mathcal{W}_w&:=\{(t,w,m,y) \in \mathcal{O}': (t,z^\ast(t,w,m,y),m,y) \in \mathcal{W}_z \},\\
\mathcal{I}_w &:=\{(t,w,m,y) \in  \mathcal{O}': (t,z^\ast(t,w,m,y),m,y) \in \mathcal{I}_z \}.
\end{aligned}
\right.
\end{equation}

Then, by Lemma \ref{b} we have 
\begin{align}\label{neww}
\begin{cases}
	\mathcal{W}_w \cap \{t<T\} =\{(t,w,m,y) \in \mathcal{O}': -g_t<w <\widehat{b}(t,m,y) \}, \\ \mathcal{I}_w \cap \{t<T\} =\{(t,w,m,y) \in  \mathcal{O}':  w\geq \widehat{b}(t,m,y) \}
	\end{cases}
\end{align}
so that we can express the optimal retirement time in terms of the initial coordinates as: {
\refstepcounter{equation}
\begin{align*}
\tau^\ast(t,w,m,y)=\inf\{s\in[t,T):W_s^\ast\geq\widehat b(s,M_s,Y_s)\}\wedge T.\tag{\theequation}\label{optimalretirement}
\end{align*}
}
\blue{At time $T$, retirement is mandatory, so the entire admissible state space is a stopping region.}

Thanks to (\ref{newb}) and Theorem \ref{verification} we can finally express the optimal retirement threshold $\widehat{b}$ and the optimal portfolio $\pi$ in terms of $b$ and $z^\ast$, respectively.

\begin{pro}\label{integralprimal}
  \blue{For every $t<T$, $m>0$, and $y>0$,}
  \begin{align}\label{bhat}
  	\widehat{b}(t,m,y)=- y Q(b(t,m),m)-\frac{b(t,m)y}{1-\gamma} Q_z(b(t,m),m).  
  	\end{align}
  %and 
 % \begin{align*}
  %	\pi^\ast(t,w,m,y) = \frac{\theta}{\sigma } z^\ast[\widetilde{J}_x(t,x^\ast,m)\frac{x^{\ast}y}{z^\ast}\frac{2-\gamma}{(1-\gamma)^2}+\widetilde{J}_{xx}(t,x^\ast,m)\frac{(x^\ast)^2y}{z^\ast}\frac{1}{(1-\gamma)^2}]\\-\sigma_y y [q_t+[\widetilde{J}(t,x^\ast,m)+\widetilde{J}_x(t,x^\ast,m) x^\ast\frac{\gamma-\gamma^2+1}{(1-\gamma)^2}+(x^\ast)^2\frac{\gamma}{(1-\gamma)^2}\widetilde{J}_{xx}(t,x^\ast,m)]]  \end{align*}
%with $x^\ast:=(z^\ast(t,w,m,y))^\frac{1}{1-\gamma}y^{\frac{\gamma}{1-\gamma}}$.
\end{pro}
\begin{proof}
We know that $\widehat{b}(t,m,y)= w^\ast(t,b(t,m)^{1-\gamma}y^{-\gamma},m,y),$ where $w^\ast(t, \cdot, m,y)$ is the inverse function of $z^\ast(t, \cdot, m,y).$ Since $J_z(t, z^\ast(t,w,m,y),m,y) =-w-g_t$, by taking $w =w^\ast(t,z,m,y)$, computations show that
\begin{align*}
J_z(t,z,m,y) =  J_z(t,z^\ast(t,w^\ast(t,z,m,y),m,y),m,y) = -w^\ast(t,z,m,y)-g_t.
\end{align*}
Hence, from (\ref{newb}), together with (\ref{4-3})-(\ref{4-4}) we have 
\begin{align*}
\widehat{b}(t,m,y) &=w^\ast(t,b(t,m)^{1-\gamma}y^{-\gamma},m,y)= -  J_z(t,b(t,m)^{1-\gamma}y^{-\gamma},m,y)-g_t \\
&= -y\widetilde{J}(t,b(t,m),m)- \frac{b(t,m)y}{1-\gamma}\widetilde{J}_x(t,b(t,m),m)-g_t \\
&=     -y\Big(Q(b(t,m),m)-q_t\Big)-\frac{b(t,m)y}{1-\gamma}Q_z(b(t,m),m)-g_t \\
&= - y Q(b(t,m),m)-\frac{b(t,m)y}{1-\gamma}Q_z(b(t,m),m).
\end{align*}

%To prove the second statement, we notice that (in the a.e. sense) $V_w(t,w,m,y)= z^\ast, V_{ww}(t,w,m,y)=-\frac{1}{J_{zz}(t,z^\ast,m,y)}, V_{wy}=\frac{-q_t-J_{zy}(t,z^\ast,m,y)}{J_{zz}(t,z^\ast,m,y)}$ (cf.\ (\ref{deJ}) and (\ref{5-2})), which then yield 
%\begin{align*}
%&\pi^\ast(t,w,m,y) = \frac{-\theta  V_w(t,w,m,y)-\sigma_y yV_{wy}(t,w,m,y)}{\sigma V_{ww}(t,w,m,y)} \\
%&= \frac{\theta}{\sigma } z^{\ast}J_{zz}(t,z^\ast,m,y)-\sigma_y y \Big(q_t+J_{zy}(t,z^\ast,m,y)\Big) \\
%&= \frac{\theta}{\sigma } z^\ast\Big[\widetilde{J}_x(t,x^\ast,m)\frac{x^{\ast}y}{z^\ast}\frac{2-\gamma}{(1-\gamma)^2}+\widetilde{J}_{xx}(t,x^\ast,m)\frac{(x^\ast)^2y}{z^\ast}\frac{1}{(1-\gamma)^2}\Big]\\&-\sigma_y y \bigg[q_t+\Big(\widetilde{J}(t,x^\ast,m)+\widetilde{J}_x(t,x^\ast,m) x^\ast\frac{\gamma-\gamma^2+1}{(1-\gamma)^2}+(x^\ast)^2\frac{\gamma}{(1-\gamma)^2}\widetilde{J}_{xx}(t,x^\ast,m)\Big)\bigg],\ a.e. \ \text{on} \ \mathcal{O}',\end{align*}
%where $x^\ast:=x^\ast(t,w,m,y)=(z^\ast(t,w,m,y))^\frac{1}{1-\gamma}y^{\frac{\gamma}{1-\gamma}}$.
\end{proof}

\begin{remark}
	Given $\widehat{b}$ as in Proposition \ref{integralprimal}, (\ref{optimalretirement}) rewrites as {$\tau^\ast(t,w,m,y)=\inf\{s\in[t,T):W_s^\ast/Y_s\geq\Upsilon(s,M_s)\}\wedge T$}, where, for $t<T$, $\Upsilon(t,m):=-Q(b(t,m),m)-\frac{b(t,m)}{1-\gamma}Q_z(b(t,m),m)$. Thus, there exists a critical wealth-to-wage ratio above which it is optimal to retire, consistently with \citet{dybvig2010lifetime}.
\end{remark}

%\section{Unhedgeable income risk}

%\citet{bensoussan2016unemployment}, \citet{jang2020optimal}

\section{Numerical Study}\label{sec6}

In this section, we present numerical illustrations of the optimal strategies. We first show the retirement boundary in the primal variables and then investigate its sensitivity to the model parameters. We use the recursive integration method proposed by \citet{huang1996pricing} to solve the integral equation (\ref{integralb2}); Appendix \ref{numerical} describes the method in detail. \blue{The numerical calculations were performed with Mathematica 13.3. The code is available  upon request .}

\begin{table}[htbp]
\caption{Basic parameters set in the numerical illustrations }
\label{tab1}
\centering
	\begin{tabular}{cccccccccc}
		\toprule  %添加表格头部粗线
		%\multicolumn{8}{c}{ }	\\
		 $\mu$	& $\sigma$ &$r$ & $\rho$   &$\gamma$ & 
		  $\mu_y$ & $\sigma_y$ &$a$ & $T$ & $K$ \\
		  \midrule
	0.08	&0.2 &0.04 &0.01& 3 &0.01&0.05 & 1/10.5 &10&2\\
	
		\bottomrule
		%\midrule  %添加表格中横线
	\end{tabular}
\end{table}

The basic parameters are listed in Table \ref{tab1}{\footnote{We consider $\gamma>1$ in the numerical study; the case $\gamma<1$ can be treated by the same method.}}. For an individual aged $x_0$, we use the Gompertz specification in \citet{milevsky2007annuitization}: $M_s=a_0^{-1} e^{(x_0+s-t-m_0)/a_0}$, with modal value $m_0=88.18$ and dispersion parameter $a_0=1/a=10.5$. The agent is 55 years old at the initial time ($x_0=55$), and the mandatory retirement age is 65 ($T=10$); hence $m=a_0^{-1}e^{(x_0-m_0)/a_0}\approx0.004.$ 

\subsection{\texorpdfstring{\blue{Optimal retirement boundaries and sensitivity analysis}}{Optimal retirement boundaries and sensitivity analysis}} 
 
\blue{We now fix the current mortality state at $m=0.004$ and illustrate the primal retirement boundary $\widehat{b}(t,m,y)$ in Figure \ref{p3}; the stopping region $\mathcal{I}_w$ defined in (\ref{neww}) lies above the surface.} The surface is linear in the income direction because the homothetic model makes the critical wealth level proportional to $y$ (cf.\ Proposition \ref{integralprimal}). 
 
   \begin{figure}[htbp]
		\setlength{\abovecaptionskip}{0pt}
	\setlength{\belowcaptionskip}{5pt}
	\begin{minipage}[b]{0.5\linewidth}
		\centering
		\includegraphics[width=2.7in]{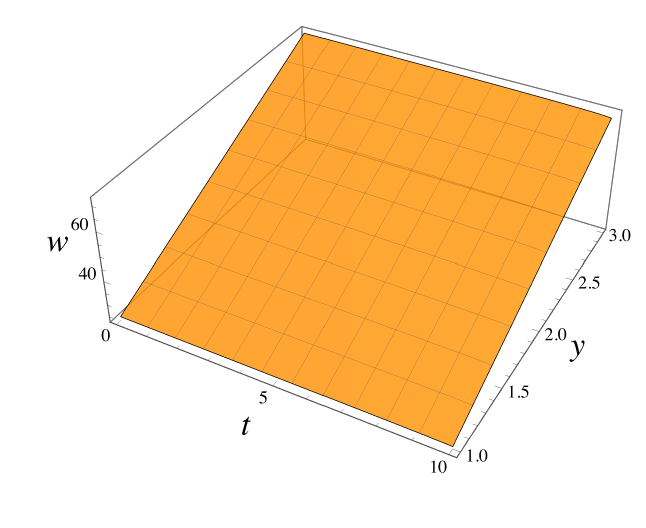}
		\caption{\blue{The surface $w=\widehat{b}(t,m,y)$ over $(t,y)$ at $m=0.004$.}}
		\label{p3}
	\end{minipage}%
	\begin{minipage}[b]{0.5\linewidth}
		\centering
		\includegraphics[width=2.7in]{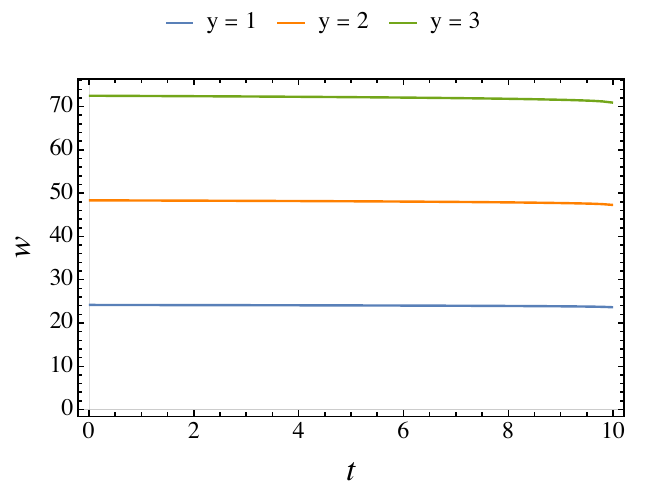}
		\caption{\blue{The primal retirement boundary at $m=0.004$ for $y=1,2,3$.}}
		\label{p4}
		\end{minipage}
		\end{figure} 

\blue{Figure \ref{p4} isolates the income channel. Doubling or tripling current income doubles or triples the boundary exactly. Economically, a larger $y$ raises human capital, namely the present value of future wages in (\ref{b(s)}), so a larger stock of financial wealth is required to compensate the agent for giving up work. Equivalently, at fixed wealth a higher income lowers the wealth-to-income ratio and makes continued employment more valuable.}

  \begin{figure}[htbp]
		\setlength{\abovecaptionskip}{0pt}
	\setlength{\belowcaptionskip}{5pt}
	
	\begin{minipage}[b]{0.5\linewidth}
		\centering
		\includegraphics[width=2.7in]{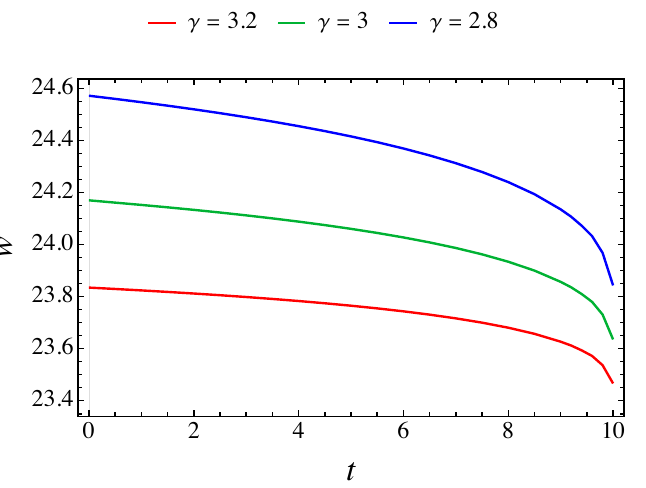}
		\caption{\blue{The primal retirement boundary at $m=0.004$ and $y=1$ for $\gamma=2.8,3,3.2$.}}
		\label{p5}
	\end{minipage}%
	\begin{minipage}[b]{0.5\linewidth}
		\centering
		\includegraphics[width=2.7in]{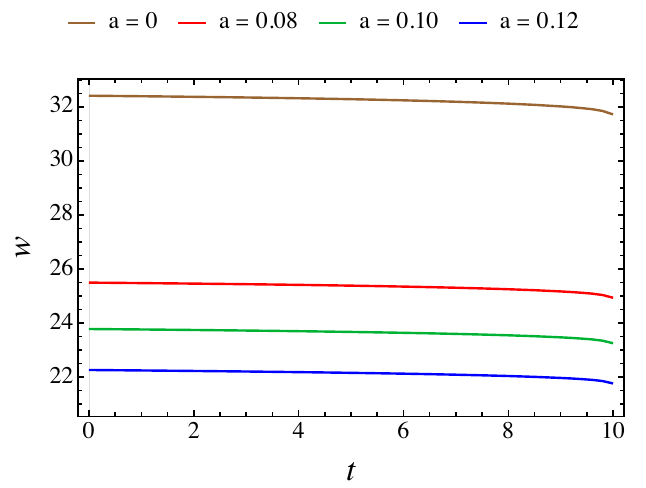}
			\caption{\blue{The primal retirement boundary at $m=0.004$ and $y=1$ for $a=0,0.08,0.10,0.12$. }}
		\label{p6}
	\end{minipage}
\end{figure}
\blue{Figure \ref{p5} shows that a larger risk-aversion parameter $\gamma$ lowers the wealth threshold. Since the stopping region lies above the boundary, an agent at a fixed state is then more likely to retire. Figure \ref{p6} isolates the Gompertz growth parameter while holding current mortality fixed. A larger $a$ makes future mortality rise more rapidly and shortens the survival-weighted utility horizon. It therefore reduces the utility contribution of consumption financed by future work and thereby lowers the retirement threshold.}

\begin{figure}[htbp]
		\setlength{\abovecaptionskip}{0pt}
	\setlength{\belowcaptionskip}{5pt}
	\begin{minipage}[b]{0.5\linewidth}
		\centering
		\includegraphics[width=2.7in]{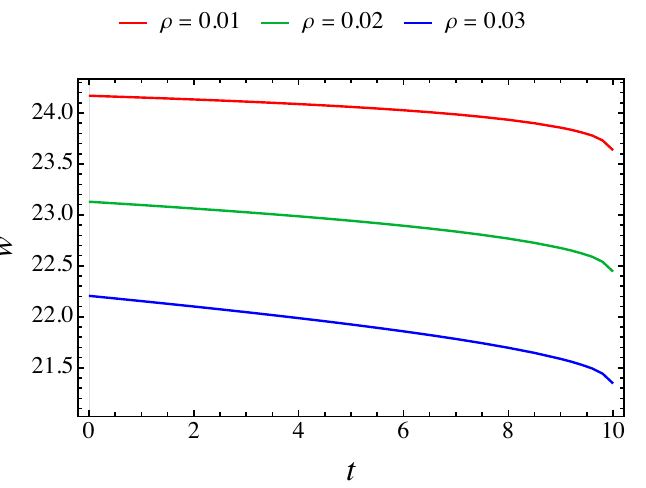}
		\caption{\blue{The primal retirement boundary at $m=0.004$ and $y=1$ for $\rho=0.01,0.02,0.03$.}}
		\label{p8}
	\end{minipage}%
	\begin{minipage}[b]{0.5\linewidth}
		\centering
		\includegraphics[width=2.7in]{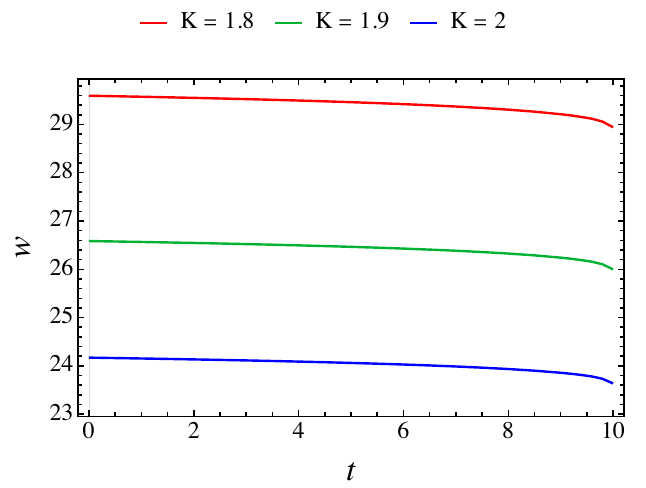}
		\caption{\blue{The primal retirement boundary at $m=0.004$ and $y=1$ for $K=1.8,1.9,2$.}}
		\label{p9}
	\end{minipage}
\end{figure}

\blue{Figure \ref{p8} shows that a larger subjective discount rate $\rho$ lowers the boundary. Greater impatience reduces the subjective weight placed on future consumption and hence the utility gain from remaining employed. Figure \ref{p9} varies the post-retirement utility weight $K$. A larger $K$ makes the retired regime more attractive and lowers the amount of wealth required to trigger retirement.}

\blue{We finally report two calibrated comparisons of the investment environment. In Figure \ref{p12}, we vary the market price of risk $\theta$ while holding $r=0.04$ and $\sigma=0.20$ fixed; this is implemented by setting $\mu=0.06,0.08,0.10$. In Figure \ref{p13}, we instead vary the constant interest rate while holding $\mu=0.08$ fixed, so $\theta=(\mu-r)/\sigma$ changes endogenously.}

\begin{figure}[htbp]
	\setlength{\abovecaptionskip}{0pt}
	\setlength{\belowcaptionskip}{5pt}
	\begin{minipage}[b]{0.5\linewidth}
		\centering
		\includegraphics[width=2.7in]{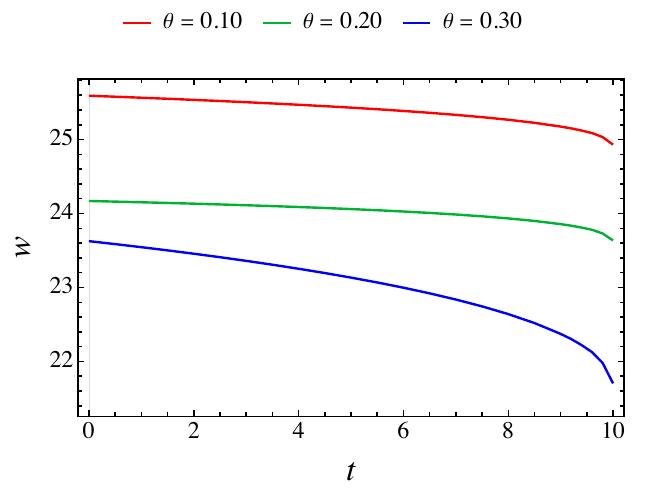}
		\caption{\blue{The primal retirement boundary at $m=0.004$ and $y=1$.}}
		\label{p12}
	\end{minipage}%
	\begin{minipage}[b]{0.5\linewidth}
		\centering
		\includegraphics[width=2.7in]{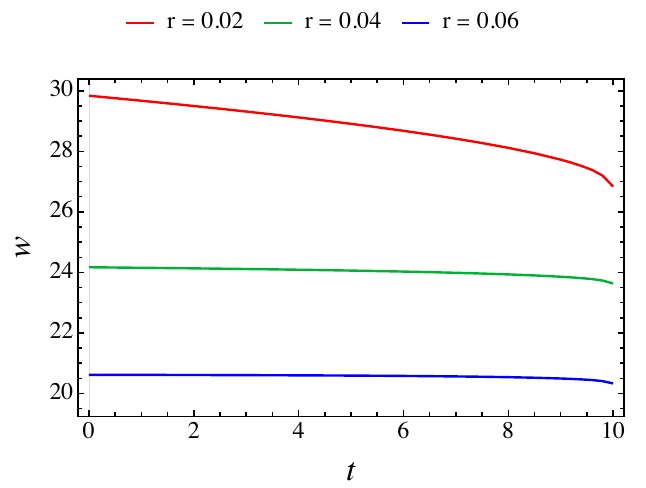}
		\caption{\blue{The primal retirement boundary at $m=0.004$ and $y=1$, with $\mu=0.08$ fixed.}}
		\label{p13}
	\end{minipage}
\end{figure}

\blue{Figure \ref{p12} shows that the boundary decreases as $\theta$ rises. A larger $\theta$ improves the compensation for financial risk, but it affects both the working and retired regimes. In addition, because labor income and the stock are exposed to the same Brownian motion, a larger $\theta$ raises the human-capital discount rate $\kappa=r-\mu_y+\sigma_y\theta$ and lowers the value of future wages. The lower boundary is therefore a joint effect, with the human-capital channel making continued work less valuable. \blue{Figure \ref{p13} also shows that the boundary decreases as $r$ rises. A higher risk-free return helps finance retirement consumption from accumulated financial wealth. Moreover, with $\mu$ fixed, $\partial_r\kappa=1-\sigma_y/\sigma=0.75>0$ under the baseline parameters, so the present value of future wages falls. Giving up these wages therefore becomes less costly, further favoring retirement at a lower wealth level. The equity premium $\mu-r$ also falls, so the investment-opportunity effect is not unambiguously favorable. The plotted decline is the net effect in this numerical example.}}

%\subsection{Comparison with constant force of mortality}
\subsection{Optimal consumption and portfolio}

\blue{To compare policy choices across wealth on a common scale, Figures \ref{p10} and \ref{p11} report the consumption--wealth ratio and the risky-asset share, respectively, at $t=0$, $m=0.004$, and $y=1$. The vertical dashed line is the computed retirement threshold $w=24.170$.}

\begin{figure}[htbp]
		\setlength{\abovecaptionskip}{0pt}
	\setlength{\belowcaptionskip}{5pt}
	\begin{minipage}[b]{0.5\linewidth}
		\centering
		\includegraphics[width=2.7in]{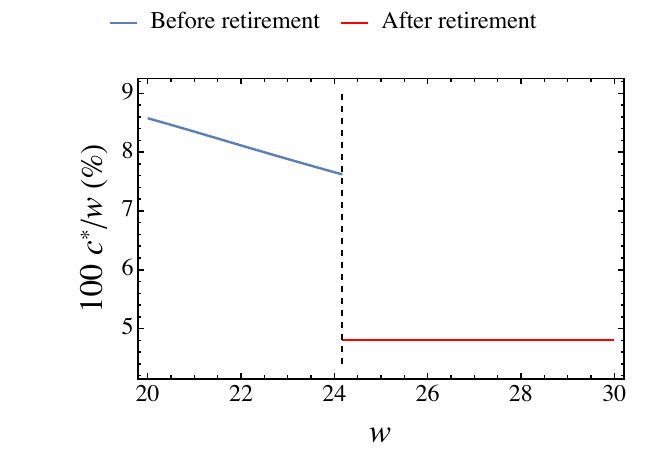}
		\caption{\blue{The optimal consumption--wealth ratio at $t=0$, $m=0.004$, and $y=1$. }}
		\label{p10}
	\end{minipage}%
	\begin{minipage}[b]{0.5\linewidth}
		\centering
		\includegraphics[width=2.7in]{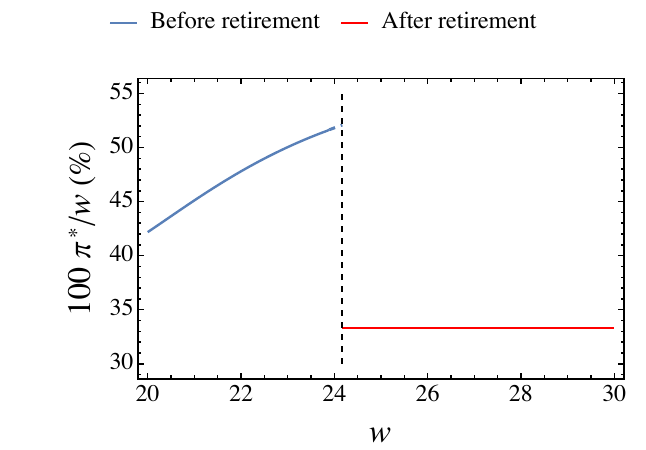}
			\caption{\blue{The optimal risky-asset share at $t=0$, $m=0.004$, and $y=1$. }}
		\label{p11}
	\end{minipage}
\end{figure}

\blue{Figure \ref{p10} shows that before retirement the consumption--wealth ratio falls from $8.578\%$ at $w=20$ to $7.622\%$ at the threshold. With labor income fixed at $y=1$, financial wealth alone does not scale the pre-retirement problem: wage income and the option to continue working support consumption when financial wealth is low, and their importance relative to wealth declines as $w$ rises. Immediately after retirement, labor income disappears. At the fixed mortality level $m=0.004$ used in the figure, the homothetic post-retirement problem gives a 
consumption-wealth ratio of $4.802\%$, independently of wealth. Thus, at the threshold, consumption falls from $1.842$ to $1.161$, equivalently a $37.0\%$ decline in the consumption--wealth ratio. This is the retirement-consumption puzzle documented by \citet{banks1998there}; as in \citet{dybvig2010lifetime}, the discrete change is generated by the switch in the marginal-utility schedule at retirement.}

\blue{Figure \ref{p11} shows the risky-asset share increasing from $42.2\%$ at $w=20$ to $51.6\%$ at $w=24$. Before retirement, the option to postpone retirement after adverse returns provides an additional labor-flexibility channel and supports greater exposure to the risky asset over the displayed wealth range, despite the positive wage--stock correlation; see \citet{dybvig2010lifetime}. At retirement, the calibrated optimal share drops from the pre-retirement boundary value of $52.08\%$ to the exact post-retirement Merton share $\theta/(\sigma\gamma)=1/3$. This discontinuity reflects a genuine regime switch: retirement removes labor income and the option to absorb adverse returns by working longer.}

\section{\texorpdfstring{\blue{Discussion and Conclusion}}{Discussion and Conclusion}}\label{sec7}

\blue{This paper characterizes optimal consumption, risky investment, and retirement timing under a Gompertz force of mortality. Duality and homogeneity reduce the control--stopping problem to an optimal stopping problem for a geometric diffusion whose drift depends on the deterministic mortality path. The dual retirement boundary is finite and strictly positive and is locally Lipschitz along each fixed Gompertz characteristic. An exact comparison across mortality characteristics supplies the additional mortality dependence used in the terminal-boundary argument. We characterize the boundary through a nonlinear integral equation and recover the optimal retirement time and admissible consumption and portfolio policies in the original variables.} \blue{Moreover, the resulting policy has a transparent economic interpretation: retirement occurs when wealth relative to labor income reaches an endogenous threshold. In the numerical example, retirement produces discrete falls in both the consumption--wealth ratio and the risky-asset share. }

\blue{The Gompertz specification is useful at several distinct stages. The initial reduction of the random lifetime requires only a deterministic, locally integrable force of mortality. The subsequent boundary analysis uses the explicit Gompertz flow, its monotonicity, and the resulting integrability estimates. Extending the full characterization to a general deterministic mortality curve would therefore require assumptions that reproduce these properties. Stochastic mortality would introduce an additional state and, unless mortality risk were spanned, an unhedgeable risk factor.}

\blue{We maintain a constant investment opportunity set. With stochastic interest rates, both the state-price density and human capital would become rate-dependent, and the interest rate would enter the stopping problem as an additional state. If rate risk were fully spanned by a sufficiently rich bond market, martingale duality might remain available, but the reduced stopping problem would generally be multidimensional. If interest-rate or labor-income risk were unspanned, the analysis would require an incomplete-market formulation, as in related retirement models such as \citet{bensoussan2016unemployment}.  We leave this extension for future research.}

\appendix
\setcounter{subsection}{0}
\renewcommand\thesubsection{A.\arabic{subsection}}
%\section{Technical Estimates}
\setcounter{equation}{0}
\renewcommand\theequation{A.\arabic{equation}}
\setcounter{lemma}{0}
    \renewcommand{\thelemma}{\Alph{section}.\arabic{lemma}}

\section{Proofs}
\renewcommand\theequation{A.\arabic{equation}}
\renewcommand\thesubsection{A.\arabic{subsection}}

\subsection{Proof of Proposition \ref{Wregularity}}
\begin{proof}\label{proW}

The convex dual of $U_2$ and the inverse of its marginal utility are
\begin{align*}
U_2^\ast(z)&=K^{\frac{1-\gamma}{\gamma}}\frac{\gamma}{1-\gamma}z^{\frac{\gamma-1}{\gamma}},\quad
\mathcal I_2^u(z)=K^{\frac{1-\gamma}{\gamma}}z^{-\frac{1}{\gamma}},\quad z>0.
\end{align*}
By (\ref{2-1}), $M_s=me^{a(s-t)}$ for $s\geq t$. Set
\begin{align*}
\Lambda:=-\frac{\rho}{\gamma}-\frac{\gamma-1}{\gamma}r-\frac{\gamma-1}{2\gamma^2}\theta^2.
\end{align*}
Using (\ref{3-6}), we obtain
\begin{align*}
\mathbb E_{t,z,m}\bigg[P_s^{\frac{\gamma-1}{\gamma}}\bigg]
&=\exp\bigg\{\frac{\gamma-1}{\gamma}\int_t^s(\rho-r+M_u)du
+\frac{1}{2}\bigg[\bigg(\frac{\gamma-1}{\gamma}\bigg)^2-\frac{\gamma-1}{\gamma}\bigg]\theta^2(s-t)\bigg\}.
\end{align*}
Since $\int_t^sM_udu=\frac{m}{a}(e^{a(s-t)}-1)$, it follows that
\begin{align*}
\mathbb E_{t,z,m}\bigg[e^{-\int_t^s(\rho+M_u)du}U_2^\ast(Z_s)\bigg]=\frac{\gamma}{1-\gamma}K^{\frac{1-\gamma}{\gamma}}z^{\frac{\gamma-1}{\gamma}}
\exp\bigg\{\Lambda(s-t)-\frac{m}{a\gamma}(e^{a(s-t)}-1)\bigg\}.
\end{align*}
Hence, Tonelli's theorem, applied to $U_2^\ast$ when $\gamma<1$ and to $-U_2^\ast$ when $\gamma>1$, and the change of variable $u=s-t$ yield
\begin{align}\label{w}
Q(t,z,m)&=\frac{\gamma}{1-\gamma}K^{\frac{1-\gamma}{\gamma}}z^{\frac{\gamma-1}{\gamma}}G(m),\qquad
G(m):=\int_0^\infty\exp\bigg\{\Lambda u-\frac{m}{a\gamma}(e^{au}-1)\bigg\}du.
\end{align}
In particular, $Q$ is independent of $t$. {\color{black}Indeed, $\lim_{u\to\infty}\frac{e^{au}-1}{u}=\infty$. Since $\frac{m}{a\gamma}>0$, we may choose $R>0$ such that
\begin{align*}
\frac{m}{a\gamma}(e^{au}-1)\geq(|\Lambda|+1)u,\quad u\geq R.
\end{align*}
On $[0,R]$, the integrand defining $G(m)$ is bounded by $e^{|\Lambda|R}$, so its integral over this interval is at most $Re^{|\Lambda|R}$. For $u\geq R$, the preceding inequality and $\Lambda\leq|\Lambda|$ give
\begin{align*}
\Lambda u-\frac{m}{a\gamma}(e^{au}-1)\leq\Lambda u-(|\Lambda|+1)u\leq-u.
\end{align*}
Hence, the integral is bounded by $\int_R^\infty e^{-u}du=e^{-R}$. This proves that $G(m)<\infty$. Therefore, $Q$ is finite.}

Finally, it is easy to see that $G\in C^1(\mathbb R_+)$ and
\begin{align*}
G'(m)=-\frac{1}{a\gamma}\int_0^\infty(e^{au}-1)\exp\bigg\{\Lambda u-\frac{m}{a\gamma}(e^{au}-1)\bigg\}du.
\end{align*}
Together with the explicit power of $z$ in (\ref{w}), this proves that $Q,Q_z\in C^{2,1}(\mathbb R^2_+)$. Moreover,
\begin{align*}
Q_{zz}(z,m)=\frac{1}{\gamma}K^{\frac{1-\gamma}{\gamma}}z^{-\frac{1}{\gamma}-1}G(m)>0,
\end{align*}
so $Q$ is strictly convex in $z$. Direct calculations prove (\ref{3-6-1}).
\end{proof}

\subsection{Proof of Theorem \ref{dualityrelation}}\label{produal}
Before we present the proof of Theorem \ref{dualityrelation}, we need the following lemmata.
\begin{lemma}\label{postbudget}
{Let $w>0$ and let $c$ be a nonnegative $\mathbb F$-progressively measurable process satisfying}
\begin{align*}
{\mathbb E_{t,w,m}\bigg[\int_t^\infty\xi_{s,t}c_sds\bigg]=w.}
\end{align*}
{Then there exists an $\mathbb F$-progressively measurable portfolio process $\pi$ such that $(c,\pi)\in\widehat{\mathcal A}(t,w,m)$. The associated wealth process is}
\begin{align}\label{postwealthbudget}
{W_s=\frac{1}{\xi_{s,t}}\mathbb E_{t,w,m}\bigg[\int_s^\infty\xi_{u,t}c_udu\bigg|\mathcal F_s\bigg],\qquad s\geq t.}
\end{align}
\end{lemma}
\begin{proof}
The proof is similar to Theorem 3.9.4 in \cite{karatzas1998methods}, and we thus omit details.
\end{proof}

\begin{lemma}\label{postmultiplier}
{For fixed $m>0$, define}
\begin{align*}
{\mathcal X(z,m):=\mathbb E_{t,z,m}\bigg[\int_t^\infty\xi_{s,t}\mathcal I_2^u(zP_s)ds\bigg],\qquad z>0.}
\end{align*}
{Then $\mathcal X(z,m)=-Q_z(z,m)$. The map $z\mapsto\mathcal X(z,m)$ is continuous and strictly decreasing, with limits $+\infty$ as $z\downarrow0$ and $0$ as $z\uparrow\infty$. Hence, for every $w>0$, there exists a unique $z^\ast=z^\ast(w,m)>0$ such that $\mathcal X(z^\ast,m)=w$.}
\end{lemma}
\begin{proof}
Since $\mathcal I_2^u(z)=K^{\frac{1-\gamma}{\gamma}}z^{-\frac{1}{\gamma}}$, the identity $e^{-\int_t^s(\rho+M_u)du}P_s=\xi_{s,t}$ and differentiation under the integral in (\ref{3-8}) give
\begin{align}\label{3-7-2}
-Q_z(z,m)&=\mathbb E_{t,z,m}\bigg[\int_t^\infty e^{-\int_t^s(\rho+M_u)du}\mathcal I_2^u(zP_s)P_sds\bigg]=\mathcal X(z,m).
\end{align}
Moreover, (\ref{w}) yields
\begin{align}\label{eqa-6-1}
\mathcal X(z,m)=-Q_z(z,m)=K^{\frac{1-\gamma}{\gamma}}z^{-\frac{1}{\gamma}}G(m).
\end{align}
The remaining claims follow directly from this expression.
\end{proof}

\begin{proof}[Proof of Theorem \ref{dualityrelation}]
\color{black}
 First, for every $(c,\pi)\in\widehat{\mathcal A}(t,w,m)$, inequality (\ref{3-5}) gives, for each $z>0$,
\begin{align}\label{eqa.3}
\widehat V(t,w,m)\leq Q(z,m)+zw.
\end{align}
Consequently,
\begin{align*}
\widehat V(t,w,m)&\leq\inf_{z>0}[Q(z,m)+zw],\qquad Q(z,m)\geq\sup_{w>0}[\widehat V(t,w,m)-zw].
\end{align*}

We next prove the reverse inequalities by constructing a strategy that attains (\ref{eqa.3}). Fix $w>0$ and let $z^\ast=z^\ast(w,m)$ be the unique multiplier in Lemma \ref{postmultiplier}. Let $Z^\ast$ solve (\ref{3-7}) with $Z_t^\ast=z^\ast$ and define
\begin{align}\label{3-7-1}
c_s^\ast:=\mathcal I_2^u(Z_s^\ast)=K^{\frac{1-\gamma}{\gamma}}(Z_s^\ast)^{-\frac{1}{\gamma}},\qquad s\geq t.
\end{align}
Because $Z_s^\ast=z^\ast P_s$, (\ref{3-7-2}) and the choice of $z^\ast$ give the exact budget equality
\begin{align*}
\mathbb E_{t,w,m}\bigg[\int_t^\infty\xi_{s,t}c_s^\ast ds\bigg]=\mathcal X(z^\ast,m)=w.
\end{align*}
Lemma \ref{postbudget} therefore provides a portfolio $\pi^\ast$ such that $(c^\ast,\pi^\ast)\in\widehat{\mathcal A}(t,w,m)$.

The equality condition in the definition of the convex dual is
\begin{align*}
U_2(c_s^\ast)-Z_s^\ast c_s^\ast=U_2^\ast(Z_s^\ast).
\end{align*}
Moreover, $e^{-\int_t^s(\rho+M_u)du}Z_s^\ast=z^\ast\xi_{s,t}$. Integrating the preceding equality and using the exact budget equality yield
\begin{align*}
\mathbb E_{t,w,m}\bigg[\int_t^\infty e^{-\int_t^s(\rho+M_u)du}U_2(c_s^\ast)ds\bigg]=Q(z^\ast,m)+z^\ast w.
\end{align*}
Thus the admissible pair $(c^\ast,\pi^\ast)$ attains the upper bound, and
\begin{align*}
\widehat V(t,w,m)=\inf_{z>0}[Q(z,m)+zw]=Q(z^\ast,m)+z^\ast w.
\end{align*}
For a fixed $z>0$, choose $w=\mathcal X(z,m)$. The same construction gives $\widehat V(t,w,m)=Q(z,m)+zw$, so the second inequality above becomes an equality:
\begin{align*}
Q(z,m)=\sup_{w>0}[\widehat V(t,w,m)-zw].
\end{align*}

It remains to identify the optimal wealth and portfolio. By (\ref{postwealthbudget}) and the strong Markov property of $(Z^\ast,M)$, for every $s\geq t$ we have
\begin{align*}
W_s^\ast=\frac{1}{\xi_{s,t}}\mathbb E\bigg[\int_s^\infty\xi_{u,t}c_u^\ast du\bigg|\mathcal F_s\bigg]=\mathbb E_{s,Z_s^\ast,M_s}\bigg[\int_s^\infty\xi_{u,s}\mathcal I_2^u(Z_u^\ast)du\bigg]=-Q_z(Z_s^\ast,M_s).
\end{align*}
By Proposition \ref{Wregularity}, $Q_z\in C^{2,1}(\mathbb R_+^2)$. Differentiating (\ref{3-6-1}) with respect to $z$, and using $U_2^{\ast\prime}(z)=-\mathcal I_2^u(z)$, gives
\begin{align*}
\frac{1}{2}\theta^2z^2Q_{zzz}+(\theta^2+\rho-r+m)zQ_{zz}+amQ_{zm}-rQ_z-\mathcal I_2^u(z)=0.
\end{align*}
It\^o's formula applied to $W_s^\ast=-Q_z(Z_s^\ast,M_s)$ then yields
\begin{align*}
dW_s^\ast=[rW_s^\ast+\pi_s^\ast(\mu-r)-c_s^\ast]ds+\sigma\pi_s^\ast dB_s,\qquad \pi_s^\ast=\frac{\theta}{\sigma}Z_s^\ast Q_{zz}(Z_s^\ast,M_s).
\end{align*}
Finally, at $w=0$, the values stated before (\ref{3-5}) and (\ref{w}) give the first dual relation directly.  This completes the proof.
\end{proof}

\subsection{Proof of Lemma \ref{reduce}}\label{proreduce}
\begin{proof}
	First, from (\ref{3-12-1}) and (\ref{3-6}) we compute the convex dual functions $U_1^\ast(z)$ and $U^\ast_2(z)$:
\begin{align}\label{4-2-1}
U^\ast_1(z)= \frac{\gamma}{1-\gamma}z^{\frac{\gamma-1}{\gamma}}, \quad U^\ast_2(z)= \frac{\gamma}{1-\gamma}K^{\frac{1-\gamma}{\gamma}}z^{\frac{\gamma-1}{\gamma}}=K^{\frac{1-\gamma}{\gamma}}U^\ast_1(z), \quad z\geq 0. 
\end{align}
{Set $\vartheta:=\frac{\gamma-1}{\gamma}$. Since $Z_s=zP_s$ and $U_2^\ast$ is homogeneous of degree $\vartheta$,}
\begin{align*}
{Q(z,m)=\mathbb E\bigg[\int_t^\infty e^{-\int_t^s(\rho+M_u)du}U_2^\ast(zP_s)ds\bigg]=z^\vartheta Q(1,m).}
\end{align*}
{For $x=z^{\frac{1}{1-\gamma}}y^{\frac{\gamma}{1-\gamma}}$, one has $x^\vartheta=z^{-\frac{1}{\gamma}}y^{-1}$ and hence $zyQ(x,m)=zyx^\vartheta Q(1,m)=z^\vartheta Q(1,m)=Q(z,m)$. Therefore,}
\begin{align}\label{4-4-1}
{Q(x,m)=\frac{Q(z,m)}{zy}.}
\end{align}

Direct calculations (cf.\ (\ref{3-7}) and (\ref{2-4-1})) show 
\begin{align}\label{4-5-1}
Z_s Y_s &= z P_s Y_s= ze^{\int^s_t (-r -\frac{1}{2}\theta^2+\rho+ M_u) du  -\int^s_t \theta dB_u} y e^{(\mu_y -\frac{1}{2}\sigma_y^2)(s-t)+ \sigma_y (B_s-B_t)}\nonumber \\
&=z y e^{\int^s_t (-r -\frac{1}{2}\theta^2+\rho+ M_u +\mu_y -\frac{1}{2}\sigma_y^2  ) du  -\int^s_t (\theta-\sigma_y) dB_u  }\nonumber \\
&=z y e^{\int^s_t (\rho+ M_u-\kappa) du} \zeta_s.
\end{align}
Then, by the continuous-time Bayes' rule (cf.\ Lemma 3.5.3 in \citet{karatzas2012brownian}), we have \blue{for any bounded stopping time $\tau\in\mathcal S$} and for $x=z^{\frac{1}{1-\gamma}}y^{\frac{\gamma}{1-\gamma}}$,
\begin{align}\label{4-6}
&\mathbb{E}_{t,z,m,y}[e^{-\int^{\tau}_t (\rho+ M_u) du}(Q( Z_{\tau}, M_{\tau})-Z_{\tau}g_\tau ) ]\nonumber\\
&=\mathbb{E}_{t,z,m,y}[e^{-\int^{\tau}_t (\rho+ M_u) du}(Q( Z_{\tau}, M_{\tau})-Z_{\tau}Y_\tau q_\tau ) ]\nonumber\\
&=\mathbb{E}_{t,z,m,y}[e^{-\int^{\tau}_t (\rho+ M_u) du}Z_\tau Y_\tau (Q( Z_{\tau}, M_{\tau})/(Z_\tau Y_\tau)- q_\tau ) ]\nonumber\\
&=\mathbb{E}_{t,z,m,y}[e^{-\int^{\tau}_t (\rho+ M_u) du}Z_\tau Y_\tau (Q( X_\tau, M_\tau)- q_\tau ) ]\nonumber \\
&=zy\mathbb{E}_{t,z,m,y}[e^{-\int^{\tau}_t (\rho+ M_u) du} e^{\int^\tau_t (\rho+ M_u-\kappa) du} \zeta_\tau (Q( X_\tau, M_\tau)- q_\tau ) ]\nonumber\\
&=zy\widehat{\mathbb{E}}_{t,x,m}[e^{-\kappa(\tau-t)}  (Q(X_\tau, M_\tau)- q_\tau ) ],
\end{align}
where we have used (\ref{4-4-1}), the definition of $X$ and (\ref{4-5-1}).
Based on the same arguments, we also have for $x=z^{\frac{1}{1-\gamma}}y^{\frac{\gamma}{1-\gamma}}$
\begin{align}\label{4-7}
\mathbb{E}_{t,z,m,y}\bigg[    \int_t^{\tau }  e^{-\int^s_t (\rho+ M_u) du} {U}^\ast_1(Z_s) ds\bigg]&= \mathbb{E}_{t,z,m,y}\bigg[    \int_t^{\tau }  e^{-\int^s_t (\rho+ M_u) du} {U}^\ast_1(X_s) Z_sY_s ds\bigg]\nonumber\\
&=  z y\mathbb{E}_{t,z,m,y}\bigg[    \int_t^{\tau }  e^{-\int^s_t (\rho+ M_u) du} {U}^\ast_1(X_s) e^{\int^s_t (\rho+ M_u-\kappa) du} \zeta_s ds\bigg]\nonumber\\
&= zy\zeta_t \widehat{\mathbb{E}}_{t,x,m}\bigg[    \int_t^{\tau }  e^{-\kappa(s-t)} {U}^\ast_1(X_s) ds\bigg]\nonumber\\
&=zy\widehat{\mathbb{E}}_{t,x,m}\bigg[    \int_t^{\tau }  e^{-\kappa(s-t)} {U}^\ast_1(X_s) ds\bigg],
\end{align}
due to $\zeta_t=1$. Here and in the sequel, $\widehat{\mathbb{E}}_{t,x,m}$ denotes the expectation under $\widehat{\mathbb{P}}$ conditioned on $X_t=x$ and $M_t=m$. Combining (\ref{4-6}) and (\ref{4-7}), together with (\ref{3-13}), we have $J(t,z,m,y)= zy \widetilde{J}(t,z^{\frac{1}{1-\gamma}}y^{\frac{\gamma}{1-\gamma}},m)$ with
\begin{align*}
\widetilde{J}(t,x,m):=\sup_{t \leq \tau \leq T
} \widehat{\mathbb{E}}_{t,x,m}\bigg[    \int_t^{\tau }  e^{-\kappa(s-t)} {U}^\ast_1(X_s) ds+ e^{-\kappa(\tau-t)}  \Big(Q(X_\tau, M_\tau)- q_\tau \Big)\bigg].
\end{align*}\end{proof}

\subsection{Proof of Proposition \ref{monotonic}}
\begin{proof}\label{promonotonic}
\textbf{Finiteness}. It is clear that $ \widehat{J}(t,x,m)\geq 0$ for all $(t,x,m) \in \mathcal{U}$. For $\gamma \in (0,1) \cup (1,\infty)$, we have $(1-K^{\frac{1-\gamma}{\gamma}})U_1^\ast(X_s)<0$, so that \begin{align*}
0&\leq \sup_{0 \leq \tau \leq T-t} \widehat{\mathbb{E}}\bigg[    \int_0^{\tau }  e^{-\kappa s}\Big( (1-K^{\frac{1-\gamma}{\gamma}}){U}^\ast_1(X_s) +1 \Big)ds \bigg]
\leq   \widehat{\mathbb{E}}\bigg[    \int_0^{T-t}  e^{-\kappa s}ds \bigg] \leq \frac{1}{\kappa}\Big(1-e^{-\kappa(T-t)}\Big).
\end{align*}

\textbf{Monotonicity in $x$}. \blue{The explicit solution (\ref{xdetailed}) gives $X_s^x=xX_s^1$. Hence the only $x$-dependent term in the expected running payoff is proportional to $(xX_s^1)^{\frac{\gamma-1}{\gamma}}$.} When $\gamma<1$, we find that $(1-K^{\frac{1-\gamma}{\gamma}}) \frac{\gamma}{1-\gamma}<0$ and $\frac{\gamma-1}{\gamma}<0$. Therefore, $x \mapsto \widehat{J}(t,x,m) $ is non-decreasing for all $(t,m)  \in [0,T] \times \mathbb{R}_+$. Analogously, when $\gamma>1$, $(1-K^{\frac{1-\gamma}{\gamma}}) \frac{\gamma}{1-\gamma}<0$ and $\frac{\gamma-1}{\gamma}>0$. Therefore,  $x \mapsto \widehat{J}(t,x,m) $ is non-increasing for all $(t,m)  \in [0,T] \times \mathbb{R}_+$.

\textbf{Monotonicity in $t$.} Monotonicity in $t$ is a trivial consequence of the fact that the $  \int_0^{\tau }  e^{-\kappa s}[ (1-K^{\frac{1-\gamma}{\gamma}}){U}^\ast_1(X_s) +1 ]ds$  is independent of $t$.

\textbf{Monotonicity in $m$.} When $\gamma<1$, we find that $\beta+1 = \frac{1}{1-\gamma}>0$, thus $m \mapsto X_t$ is increasing $\widehat{\mathbb{P}}$-a.s. for all $t \in [0,T]$ by uniqueness of the trajectories {in (\ref{xdetailed})}. As a consequence, $m \mapsto \widehat{J}(t,x,m) $ is non-decreasing for all $(t,x)  \in [0,T] \times \mathbb{R}_+.$

When $\gamma>1$, we find that $\beta+1 = \frac{1}{1-\gamma}<0$, thus $m \mapsto X_t$ is decreasing $\widehat{\mathbb{P}}$-a.s. for all $t \in [0,T]$ by uniqueness of the trajectories in (\ref{xdetailed}), and $m \mapsto \widehat{J}(t,x,m) $ is therefore non-decreasing for all $(t,x)  \in [0,T] \times \mathbb{R}_+.$

\end{proof}

\subsection{Proof of Proposition \ref{lip}}
\begin{proof}\label{prolip}
\color{black}
The proof is divided into two steps. 

Set $\vartheta:=\frac{\gamma-1}{\gamma}$. By (\ref{XN}), the payoff associated with any stopping time $\tau\leq T-t$ is
\begin{align}\label{payofflip}
\widehat{\mathbb E}_{x,m}\bigg[\int_0^\tau e^{-\kappa s}\bigg((1-K^{\frac{1-\gamma}{\gamma}})\frac{\gamma}{1-\gamma}x^\vartheta H_sN(s,m)+1\bigg)ds\bigg].
\end{align}
The process $H$ is a positive martingale with $\widehat{\mathbb E}[H_s]=1$. Moreover, on $[0,T]$ the functions $N$ and $N_m$ satisfy
\begin{align}\label{Nm}
N_m(s,m)=-\frac{1}{\gamma}\frac{e^{as}-1}{a}N(s,m),\qquad \frac{N(s,m+\epsilon)-N(s,m)}{\epsilon}=\int_0^1N_m(s,m+\lambda\epsilon)d\lambda,
\end{align}
and $N(s,m)\leq C$ uniformly over $s\in[0,T]$ and $m>0$.

\textbf{Step 1: the $x$- and $m$-variables.} The $x$-variable is simpler and can be handled by the same comparison as in the first part of Step 2 and Step 3 of Proposition 3.4 in \citet{de2019free}. Indeed, $X_s^x=xX_s^1$, and (\ref{payofflip}) gives, for $x_1,x_2$ in a compact interval contained in $\mathbb R_+$,
\begin{align*}
|\widehat J(t,x_2,m)-\widehat J(t,x_1,m)|&\leq C|x_2^\vartheta-x_1^\vartheta|\int_0^{T-t}e^{-\kappa s}N(s,m)\widehat{\mathbb E}[H_s]ds\leq C|x_2-x_1|,
\end{align*}
where the last constant is uniform on compact subsets of $\mathcal U$. The same comparison with the optimal stopping time at $(t,x,m)$, applied to $x+\epsilon$ and $x-\epsilon$, gives (\ref{lipx}) at every differentiability point. We omit the calculation.

We now consider the additional $m$-variable. First, (\ref{Nm}), the mean value theorem, and $|\sup_\tau f_\tau-\sup_\tau g_\tau|\leq\sup_\tau|f_\tau-g_\tau|$ give, for $m_1,m_2$ in a compact interval,
\begin{align*}
|\widehat J(t,x,m_2)-\widehat J(t,x,m_1)|&\leq Cx^\vartheta|m_2-m_1|\int_0^{T-t}e^{-\kappa s}\frac{e^{as}-1}{a}ds\leq C|m_2-m_1|.
\end{align*}
Hence $\widehat J(t,x,\cdot)$ is locally Lipschitz. Fix $(t,x,m)$ and let $\tau^\ast:=\tau^\ast(t,x,m)$. Using $\tau^\ast$ as a suboptimal rule at $m+\epsilon$ gives
\begin{align*}
\widehat J(t,x,m+\epsilon)-\widehat J(t,x,m)&\geq(1-K^{\frac{1-\gamma}{\gamma}})\frac{\gamma}{1-\gamma}x^\vartheta\widehat{\mathbb E}_{x,m}\bigg[\int_0^{\tau^\ast}e^{-\kappa s}H_s\big(N(s,m+\epsilon)-N(s,m)\big)ds\bigg].
\end{align*}
Using the same rule at $m-\epsilon$ gives
\begin{align*}
\widehat J(t,x,m)-\widehat J(t,x,m-\epsilon)&\leq(1-K^{\frac{1-\gamma}{\gamma}})\frac{\gamma}{1-\gamma}x^\vartheta\widehat{\mathbb E}_{x,m}\bigg[\int_0^{\tau^\ast}e^{-\kappa s}H_s\big(N(s,m)-N(s,m-\epsilon)\big)ds\bigg].
\end{align*}
For $m$ in a compact interval, (\ref{Nm}) gives a common integrable bound for both difference quotients. Dividing by $\epsilon$, using dominated convergence, and arguing as in Step 3 of Proposition 3.4 in \citet{de2019free}, we obtain at every differentiability point
\begin{align*}
\widehat J_m(t,x,m)&=(K^{\frac{1-\gamma}{\gamma}}-1)\frac{x^\vartheta}{1-\gamma}\widehat{\mathbb E}_{x,m}\bigg[\int_0^{\tau^\ast}e^{-\kappa s}\frac{e^{as}-1}{a}H_sN(s,m)ds\bigg],
\end{align*}
which is (\ref{lipm}).

\textbf{Step 2: the time variable.}  Set $h:=T-t$ and $\nu_\epsilon:=\tau^\ast\wedge(h-\epsilon)$ for $0<\epsilon<h$. Since $\nu_\epsilon$ is admissible at $t+\epsilon$, monotonicity gives
\begin{align*}
0&\geq\widehat J(t+\epsilon,x,m)-\widehat J(t,x,m)\nonumber\\
&\geq-\widehat{\mathbb E}_{x,m}\bigg[\int_{\nu_\epsilon}^{\tau^\ast}e^{-\kappa s}\bigg|(1-K^{\frac{1-\gamma}{\gamma}})U_1^\ast(X_s)+1\bigg|ds\bigg].
\end{align*}
Since the interval $[\nu_\epsilon,\tau^\ast]$ is empty unless $\tau^\ast\geq h-\epsilon$, and since $|(1-K^{\frac{1-\gamma}{\gamma}})U_1^\ast(X_s)+1|\leq CX_s^\vartheta+1$, Tonelli's theorem and (\ref{XN}) yield
\begin{align}\label{b5}
\widehat J(t+\epsilon,x,m)-\widehat J(t,x,m)&\geq-Cx^\vartheta\int_{h-\epsilon}^h\widetilde{\mathbb P}(\tau^\ast\geq s)ds-\int_{h-\epsilon}^h\widehat{\mathbb P}(\tau^\ast\geq s)ds\nonumber\\
&\geq-Cx^\vartheta\epsilon\widetilde{\mathbb P}(\tau^\ast\geq h-\epsilon)-\epsilon\widehat{\mathbb P}(\tau^\ast\geq h-\epsilon).
\end{align}
Probabilities are bounded by one, so this is a local Lipschitz estimate from the right. Applying the same estimate to the problem starting at $t-\epsilon$ gives the estimate from the left. Finally, Markov's inequality in (\ref{b5}) gives
\begin{align*}
\widehat J(t+\epsilon,x,m)-\widehat J(t,x,m)&\geq-\frac{\epsilon}{T-t-\epsilon}\bigg(Cx^\vartheta\widetilde{\mathbb E}[\tau^\ast]+\widehat{\mathbb E}[\tau^\ast]\bigg).
\end{align*}
At a differentiability point, division by $\epsilon$ and passage to the limit give (\ref{lipt}); the upper bound is the time monotonicity.

The estimates in Steps 1--2 are locally uniform in $(t,x,m)$. Combining the preceding estimates, we conclude that $\widehat J$ is locally Lipschitz continuous on $\mathcal U$ and differentiable almost everywhere. At every differentiability point, the preceding difference-quotient comparisons give (\ref{lipx})--(\ref{lipt}).
\end{proof}

\subsection{Proof of Proposition \ref{limit}}
\begin{proof}\label{prolimit}

Set $\vartheta=\frac{\gamma-1}{\gamma}$ and $c_\gamma:=\big(K^{\frac{1-\gamma}{\gamma}}-1\big)\frac{\gamma}{1-\gamma}>0$. By (\ref{XN}), the running payoff in (\ref{widehat{J}}) is $1-c_\gamma x^\vartheta H_sN(s,m)$. Choosing $\tau=T-t$ and using Proposition \ref{monotonic}, we obtain
\begin{align*}
\widehat{\mathbb E}_{x,m}\bigg[\int_0^{T-t}e^{-\kappa s}\big(1-c_\gamma x^\vartheta H_sN(s,m)\big)ds\bigg]&\leq\widehat J(t,x,m)\leq\int_0^{T-t}e^{-\kappa s}ds.
\end{align*}
Since $\widehat{\mathbb E}[H_s]=1$, this is equivalent to
\begin{align}\label{limitsqueeze}
\int_0^{T-t}e^{-\kappa s}ds-c_\gamma x^\vartheta\int_0^{T-t}e^{-\kappa s}N(s,m)ds&\leq\widehat J(t,x,m)\leq\int_0^{T-t}e^{-\kappa s}ds.
\end{align}
When $\gamma>1$, $\vartheta>0$ and $x^\vartheta\downarrow0$ as $x\downarrow0$. When $0<\gamma<1$, $\vartheta<0$ and $x^\vartheta\downarrow0$ as $x\uparrow\infty$. In either case, (\ref{limitsqueeze}) gives
\begin{align*}
\lim\widehat J(t,x,m)=\int_0^{T-t}e^{-\kappa s}ds=\frac{1}{\kappa}\big(1-e^{-\kappa(T-t)}\big),
\end{align*}
with the limit taken as $x\downarrow0$ when $\gamma>1$ and as $x\uparrow\infty$ when $0<\gamma<1$.

\color{black}Next, we prove $\lim_{x\uparrow\infty}\widehat J(t,x,m)=0$ when $\gamma>1$; the proof of $\lim_{x\downarrow0}\widehat J(t,x,m)=0$ when $0<\gamma<1$ is identical. For every stopping time $\tau\leq T-t$, we have
\begin{align*}
\int_0^\tau e^{-\kappa s}\big(1-c_\gamma x^\vartheta H_sN(s,m)\big)ds&\leq\int_0^{T-t}e^{-\kappa s}\big(1-c_\gamma x^\vartheta H_sN(s,m)\big)^+ds.
\end{align*}
Taking expectations and then the supremum over $\tau$ gives
\begin{align}\label{limitpositivepart}
0\leq\widehat J(t,x,m)&\leq\widehat{\mathbb E}_{x,m}\bigg[\int_0^{T-t}e^{-\kappa s}\big(1-c_\gamma x^\vartheta H_sN(s,m)\big)^+ds\bigg].
\end{align}
For every $s>0$, one has $H_sN(s,m)>0$. Hence, as $x\uparrow\infty$ when $\gamma>1$, the positive part in (\ref{limitpositivepart}) converges to zero. It is bounded by one, and the time horizon is finite, so dominated convergence proves the desired limit. When $0<\gamma<1$ and $x\downarrow0$, one again has $x^\vartheta\uparrow\infty$, and the same argument proves the remaining limit.
\end{proof}

\subsection{Proof of Lemma \ref{b}}
\begin{proof}\label{prob}

\textbf{Step 1:} Here we prove (\ref{stopping1}). Equation (\ref{stopping2}) can be proved by similar arguments. Since $x\mapsto\widehat J(t,x,m)$ is non-decreasing by Proposition \ref{monotonic}, we can define
\begin{align*}
b(t,m):=\sup\{x>0:\widehat J(t,x,m)\leq0\},\qquad\sup\emptyset:=0,
\end{align*}
so that $\mathcal I \cap \{t<T\} =\{(t,x,m)\in[0,T)\times\mathbb R_+^2:0<x\leq b(t,m)\}.$

\textbf{Step 2:} Here we prove property (i) for the case $\gamma<1$. Analogous considerations allow us to prove it when $\gamma>1$. We know that $t\mapsto\widehat J(t,x,m)$ is non-increasing by Proposition \ref{monotonic}. Therefore, for $t\geq0$ and $s\leq T-t$, we have
\begin{align*}
b(t,m)=\sup\{x>0:\widehat J(t,x,m)\leq0\}\leq\sup\{x>0:\widehat J(t+s,x,m)\leq0\}=b(t+s,m).
\end{align*}

\textbf{Step 3:} Here we prove property (ii), again only in the case $\gamma<1$. Since $m\mapsto\widehat J(t,x,m)$ is non-decreasing by Proposition \ref{monotonic}, for $m>0$ and $h>0$, we have
\begin{align*}
b(t,m)=\sup\{x>0:\widehat J(t,x,m)\leq0\}\geq\sup\{x>0:\widehat J(t,x,m+h)\leq0\}=b(t,m+h),
\end{align*}
which implies that $m\mapsto b(t,m)$ is non-increasing when $\gamma<1$.

\textbf{Step 4:} Here we provide the proof of property (iii), i.e., the bounds of the boundary. Noticing that, due to (\ref{widehat{J}}),
\begin{align*}
\mathcal R:=\{(t,x,m)\in[0,T)\times\mathbb R_+^2:(1-K^{\frac{1-\gamma}{\gamma}})U_1^\ast(x)+1>0\}\subseteq\mathcal W,
\end{align*}
we have
\begin{align}\label{stoppingbound}
\mathcal R^C:=\{(t,x,m)\in[0,T)\times\mathbb R_+^2:(1-K^{\frac{1-\gamma}{\gamma}})U_1^\ast(x)+1\leq0\}\supseteq\mathcal I.
\end{align}
We first show
\begin{align*}
b(t,m)\leq L:=\Big[\Big(K^{\frac{1-\gamma}{\gamma}}-1\Big)\frac{\gamma}{1-\gamma}\Big]^{\frac{\gamma}{1-\gamma}}
\end{align*}
for $\gamma<1$. Recalling $U_1^\ast(x)$ in (\ref{4-2-1}), one has
\begin{align*}
(1-K^{\frac{1-\gamma}{\gamma}})U_1^\ast(x)+1\leq0\quad\Longleftrightarrow\quad x^{\frac{\gamma-1}{\gamma}}\geq\frac{1}{\big(K^{\frac{1-\gamma}{\gamma}}-1\big)\frac{\gamma}{1-\gamma}},
\end{align*}
for all $(t,m)\in[0,T)\times\mathbb R_+$; that is,
\begin{align*}
(t,x,m)\in\mathcal R^C\quad\Longleftrightarrow\quad x\leq\Big[\Big(K^{\frac{1-\gamma}{\gamma}}-1\Big)\frac{\gamma}{1-\gamma}\Big]^{\frac{\gamma}{1-\gamma}}.
\end{align*}
Because (\ref{stopping1}) holds, (\ref{stoppingbound}) and the latter equivalence imply
\begin{align*}
b(t,m)\leq\Big[\Big(K^{\frac{1-\gamma}{\gamma}}-1\Big)\frac{\gamma}{1-\gamma}\Big]^{\frac{\gamma}{1-\gamma}},\qquad(t,m)\in[0,T)\times\mathbb R_+.
\end{align*}
Analogous arguments, now using $\gamma>1$, allow us to obtain
{
\begin{align*}
b(t,m)&\geq\Big[\Big(K^{\frac{1-\gamma}{\gamma}}-1\Big)\frac{\gamma}{1-\gamma}\Big]^{\frac{\gamma}{1-\gamma}},\qquad(t,m)\in[0,T)\times\mathbb R_+.\qedhere
\end{align*}
}
\end{proof}

\subsection{Proof of Theorem \ref{lipschitz} }\label{prolipschitz}

\begin{proof}
\color{black}
We apply the one-dimensional Theorem 4.3 of \citet{de2019lipschitz}. Fix $m>0$, so that $\overline m(t)=me^{at}$. Using $L$ from Lemma \ref{b}, introduce the coordinate
\begin{align*}
\widetilde z:=-\frac{\gamma-1}{\gamma}\log(x/L)-\frac{\overline m(t)}{a\gamma}.
\end{align*}
The corresponding process $\widetilde Z_s=-\frac{\gamma-1}{\gamma}\log(X_s/L)-\frac{\overline m(t+s)}{a\gamma}$ satisfies, by It\^o's formula,
\begin{align*}
d\widetilde Z_s=-\frac{\gamma-1}{\gamma}\bigg((\beta+1)(\rho-r)+\mu_1-\frac{\sigma_1^2}{2}\bigg)ds-\frac{\gamma-1}{\gamma}\sigma_1d\widehat B_s,\qquad\widetilde Z_0=\widetilde z.
\end{align*}
Both coefficients are constant, and the diffusion coefficient is nonzero by Assumption \ref{nondegenerate}. Since the running payoff in (\ref{mathcalJ}) equals $1-(x/L)^{\frac{\gamma-1}{\gamma}}$, the transformed problem is
\begin{align*}
e^{-\kappa t}\mathcal J(t,x)&=\sup_{0\leq\tau\leq T-t}\widehat{\mathbb E}_{\widetilde z}\bigg[\int_0^\tau h(t+s,\widetilde Z_s)ds\bigg],\\
h(t,\widetilde z)&:=e^{-\kappa t}\Big(1-e^{-\widetilde z-\frac{\overline m(t)}{a\gamma}}\Big).
\end{align*}
For both ranges of $\gamma$, Lemma \ref{b} identifies its stopping set by
\begin{align*}
\widetilde z\leq-\frac{\gamma-1}{\gamma}\log(\overline b(t)/L)-\frac{\overline m(t)}{a\gamma}.
\end{align*}

We verify the assumptions of Theorem 4.3 of \cite{de2019lipschitz}. In its notation, the stopping and terminal rewards are zero, so (A1), (A2), and part (i) of Corollary 3.3 hold, and the generator term associated with the stopping reward vanishes. The drift derivative is zero. Moreover,
\begin{align*}
\partial_{\widetilde z}h(t,\widetilde z)&=e^{-\kappa t-\widetilde z-\frac{\overline m(t)}{a\gamma}}>0,\\
\partial_t h(t,\widetilde z)&=-\kappa e^{-\kappa t}+\Big(\kappa+\frac{\overline m(t)}{\gamma}\Big)\partial_{\widetilde z}h(t,\widetilde z).
\end{align*}
Thus (C) and (D) hold, with a constant uniform over $t\in[0,T]$. The sign-change curve is $-\overline m(t)/(a\gamma)$. For any compact interval $I\subset[0,T)$ and any $\widetilde z_0$ above this curve on $I$,
\begin{align*}
\inf_{t\in I,\widetilde z\leq\widetilde z_0}\partial_{\widetilde z}h(t,\widetilde z)>0.
\end{align*}
This verifies parts (i)--(ii) of Theorem 4.3. The regularity and uniform-integrability requirements in Assumption 2.2 follow from finite-horizon Gaussian exponential moments, and the transformed value is continuous by Proposition \ref{lip}. 

Theorem 4.3 of \cite{de2019lipschitz} also makes the transformed threshold finite and locally Lipschitz on $(0,T)$. To cover $t=0$, extend the same formulas for $\overline m$ and $h$ to $[-\delta,T]$ and shift the time origin. The original time zero is then an interior point, while the problem at nonnegative original times is unchanged. This gives the one-sided Lipschitz estimate at zero. The inverse coordinate transformation now implies that $\overline b$ is strictly positive, finite, and locally Lipschitz on $[0,T)$. Since every current mortality state lies on one of these paths, $0<b(t,m)<\infty$ for every $t<T$ and $m>0$.
\end{proof}

\subsection{Proof of Lemma \ref{interior} }\label{prointerior}
\begin{proof}
{\color{black}
The inequality $\widehat\tau(t,x)\geq\overline\tau^\ast(t,x)$ follows directly from the definitions. We first prove immediate strict entry from a boundary point. Fix $t_0<T$ and $x_0=\overline b(t_0)$. Since $\overline b$ is locally Lipschitz and strictly positive, there are constants $C>0$ and $s_0>0$ such that
\begin{align*}
|\log\overline b(t_0+s)-\log x_0|\leq Cs,\qquad0\leq s\leq s_0.
\end{align*}
Starting from $X_0=x_0$, one has
\begin{align*}
\log\bigg(\frac{X_s}{x_0}\bigg)=\sigma_1\widehat B_s+\int_0^s\bigg[(\beta+1)(\rho-r+\overline m(t_0+u))+\mu_1-\frac{\sigma_1^2}{2}\bigg]du.
\end{align*}
The drift in the integral is bounded on $[0,s_0]$. Hence, with
\begin{align*}
D_s:=\log X_s-\log\overline b(t_0+s),
\end{align*}
there is a constant $C_0>0$ such that
\begin{align*}
D_s=\sigma_1\widehat B_s+R_s,\qquad|R_s|\leq C_0s,\qquad0\leq s\leq s_0.
\end{align*}
Assumption \ref{nondegenerate} gives $\sigma_1\neq0$. The law of the iterated logarithm, as stated in Remark 3.4 of \citet{baldi2017stochastic}, yields sequences $s_n\downarrow0$ and $r_n\downarrow0$ such that
\begin{align*}
D_{s_n}>0,\qquad D_{r_n}<0
\end{align*}
for all sufficiently large $n$, a.s. When $\gamma<1$, the inequalities $D_{r_n}<0$ show entry into $\{x<\overline b(t)\}$ at arbitrarily small times. When $\gamma>1$, the inequalities $D_{s_n}>0$ show entry into $\{x>\overline b(t)\}$ at arbitrarily small times. Thus, in both cases,
\begin{align*}
\widehat\tau(t_0,\overline b(t_0))=0,\qquad\widehat{\mathbb P}\text{-a.s.}
\end{align*}

\blue{If $t=T$, both times are zero. If $(t,x)$ lies in the strict stopping region, both times are again zero; initial boundary points have already been treated. It remains to start in the continuation region.} On $\{\overline\tau^\ast(t,x)<T-t\}$, continuity of $X$ gives $X_{\overline\tau^\ast}=\overline b(t+\overline\tau^\ast)$. The strong Markov property and the immediate-entry result applied at time $t+\overline\tau^\ast$ show that the path enters the strict stopping region at times arbitrarily close to $\overline\tau^\ast$. Hence $\widehat\tau(t,x)\leq\overline\tau^\ast(t,x)$. On $\{\overline\tau^\ast(t,x)=T-t\}$ both times equal the horizon. Together with the opposite inequality, this proves the lemma.
}
\end{proof}

 \vspace{-\baselineskip}

\subsection{Proof of Theorem \ref{verification} }\label{proofsolution}
Before we present the proof of Theorem \ref{verification}, we need the following lemma.

\begin{lemma}\label{budegt}
{For $\tau\in\mathcal S$, let $w>-g_t$, let $c\geq0$ be an $\mathbb F$-progressively measurable consumption process, and let $\phi>0$ a.s.\ be $\mathcal F_\tau$-measurable. Suppose that}
\begin{align*}
\mathbb E_{t,w,m,y}\bigg[\int_t^\tau\xi_{s,t}c_sds+\xi_{\tau,t}(\phi+g_\tau)\bigg]=w+g_t.
\end{align*}
Then there exists an $\mathbb F$-progressively measurable portfolio process $\pi$, defined up to $\tau$, such that $(c,\pi)$ is admissible before retirement,
\begin{align*}
W_s^{c,\pi,\tau}\geq-g_s,\quad t\leq s<\tau,\qquad W_\tau^{c,\pi,\tau}=\phi,
\end{align*}
and the corresponding wealth process admits the conditional budget representation
\begin{align*}
W_s^{c,\pi,\tau}+g_s=\frac{1}{\xi_{s,t}}\mathbb E_{t,w,m,y}\bigg[\int_s^\tau\xi_{u,t}c_udu+\xi_{\tau,t}(\phi+g_\tau)\bigg|\mathcal F_s\bigg],\quad t\leq s\leq\tau.
\end{align*}
\end{lemma}
\begin{proof}
\color{black}
This is the finite-horizon version of the martingale-representation argument in Lemma 6.3 of \citet{karatzas2000utility}. The details are omitted. 
%Indeed, set
%\begin{align*}
%\mathcal M_s:=\mathbb E_{t,w,m,y}\bigg[\int_t^\tau\xi_{u,t}c_udu+\xi_{\tau,t}(\phi+g_\tau)\bigg|\mathcal F_s\bigg],\quad t\leq s\leq\tau.
%\end{align*}
%\blue{The terminal random variable is nonnegative and integrable, so $\mathcal M$ is a uniformly integrable martingale. By the martingale representation theorem, there is a progressively measurable process $\psi$ such that}
%\begin{align*}
%\blue{d\mathcal M_s=\psi_sdB_s,\qquad\int_t^\tau|\psi_s|^2ds<\infty,\quad\mathbb P\text{-a.s.}}
%\end{align*}
%Define the augmented wealth by
%\begin{align*}
%W_s+g_s:=\frac{1}{\xi_{s,t}}\bigg(\mathcal M_s-\int_t^s\xi_{u,t}c_udu\bigg),\quad t\leq s\leq\tau.
%\end{align*}
%The assumed budget equality gives $W_t=w$, while the definition gives $W_\tau=\phi$. The conditional-expectation form also yields $W_s+g_s\geq0$. Recall that
%\begin{align*}
%dg_s=[(r+\theta\sigma_y)g_s-Y_s]ds+\sigma_yg_sdB_s.
%\end{align*}
%It\^o's formula for the displayed definition of $W_s+g_s$ shows that the choice
%\begin{align*}
%\pi_s:=\frac{1}{\sigma}\bigg(\frac{\psi_s}{\xi_{s,t}}+\theta(W_s+g_s)-\sigma_yg_s\bigg)
%\end{align*}
%makes $W$ satisfy (\ref{sde}) on $[t,\tau]$. On each localized finite interval, continuity of $\xi_{s,t}^{-1}$, $W_s+g_s$, and $g_s$, together with the square integrability of $\psi$, implies the required local square integrability of $\pi$. Hence the replicated strategy is admissible.
\end{proof}

%\begin{lemma}
%	For $\tau \in \mathcal{S}$ and fixed $(t,m) \in [0,T] \times \mathbb{R}_+$, we define 
%\begin{align*}
%\bar{\mathcal{X}}_\tau(z):= \mathbb{E}_{t,z,m,y}\bigg[  \int_t^{\tau}\xi_{s,t}\mathcal{I}^u_1(zP_s) ds       + \xi_{\tau,t}\mathcal{I}^{Q}(Z_\tau,M_\tau)        \bigg],
%\end{align*}
%where $\mathcal{I}^u_1$ is the inverse function of $U'_1(\cdot)$ and $\mathcal{I}^Q(z,m)$ is the inverse function of $Q_z$. Then, for any $w>0$, there exists $z>0$ such that $\bar{\mathcal{X}}_{\tau^\ast(z,m,y)}(z)=w$.
%\end{lemma}
%\begin{proof}
%Since the regularity of $Q$ in Proposition \ref{Wregularity}, $\mathcal{I}^Q$ is well-defined. Moreover, 
%\begin{align*}
%\lim_{z \to \infty}	\mathcal{I}^Q(z,m)= \lim_{z\to \infty} 
%\end{align*}

%	Using the monotone convergence theorem, we obtain 
%	\begin{align*}
%	\lim_{z\to 0}	\bar{\mathcal{X}}_\tau(t,z)=\infty, \quad \lim_{z \to \infty} \bar{\mathcal{X}}_\tau(t,z)=0.
%				\end{align*}
%				We claim that $\bar{\mathcal{X}}_\tau(t,z)$ is continuous in $z$. In fact, 
%\end{proof}

\begin{proof}[Proof of Theorem \ref{verification}]
\color{black}
First, we show the duality relations. Fix an arbitrary admissible strategy $(c,\pi,\tau)\in\mathcal A(t,w,m,y)$. For every $z>0$, (\ref{3-12}) gives
\begin{align*}
V(t,w,m,y)\leq J(t,z,m,y)+z(w+g_t).
\end{align*}
Consequently,
\begin{align*}
V(t,w,m,y)&\leq\inf_{z>0}[J(t,z,m,y)+z(w+g_t)],\\
J(t,z,m,y)&\geq\sup_{w>-g_t}[V(t,w,m,y)-z(w+g_t)].
\end{align*}

We now prove the reverse inequalities. By Lemma \ref{convex}, for every $(t,w,m,y)\in\mathcal O'$ there exists a unique $z^\ast=z^\ast(t,w,m,y)>0$ such that
\begin{align*}
w+g_t=-J_z(t,z^\ast,m,y).
\end{align*}
For $z>0$, let $\mathcal I_1^u(z):=z^{-\frac{1}{\gamma}}$ denote the inverse of $U_1'$, let $\tau_z^\ast$ be the first-entry optimal stopping time for (\ref{3-13}), and define
\begin{align*}
\overline{\mathcal X}(t,z,m,y):=\mathbb E_{t,z,m,y}\bigg[\int_t^{\tau_z^\ast}\xi_{s,t}\mathcal I_1^u(Z_s)ds+\xi_{\tau_z^\ast,t}\big(-Q_z(Z_{\tau_z^\ast},M_{\tau_z^\ast})+g_{\tau_z^\ast}\big)\bigg].
\end{align*}
For the given $(t,w,m,y)$, let $Z^\ast$ solve (\ref{3-7}) with $Z_t^\ast=z^\ast$, set $\tau^\ast:=\tau_{z^\ast}^\ast$, and introduce the candidate pre-retirement consumption process
\begin{align*}
c_s^\ast:=\mathcal I_1^u(Z_s^\ast)=(Z_s^\ast)^{-\frac{1}{\gamma}},\qquad t\leq s<\tau^\ast.
\end{align*}

We first claim that
\begin{align*}
\overline{\mathcal X}(t,z,m,y)=-J_z(t,z,m,y).
\end{align*}
For a fixed stopping time $\tau$, denote the objective in (\ref{3-13}) by $\overline J(t,z,m,y;\tau)$. Since $Z_s=zP_s$, differentiation under the expectation gives
\begin{align}\label{envelopederivative}
\frac{\partial}{\partial z}\overline J(t,z,m,y;\tau)
=\mathbb E_{t,z,m,y}\bigg[-\int_t^\tau\xi_{s,t}\mathcal I_1^u(Z_s)ds+\xi_{\tau,t}\big(Q_z(Z_\tau,M_\tau)-g_\tau\big)\bigg].
\end{align}
Indeed, $U_1^{\ast\prime}=-\mathcal I_1^u$ and $e^{-\int_t^s(\rho+M_u)du}P_s=\xi_{s,t}$. The power form of $U_1^\ast$ and $Q$, the finite horizon, and the lognormal moments of $(P,Y)$ justify the differentiation and provide a common dominator when $z$ ranges over a compact subset of $\mathbb R_+$.

\blue{It remains to verify that the derivative may be evaluated at the optimal stopping time. For the fixed $(t,m,y)$, use the mortality trajectory $s\mapsto me^{a(s-t)}$. The homogeneity relations (\ref{jandj}), (\ref{4-3}), (\ref{4-4}), and (\ref{wxm}) give $\tau_z^\ast=t+\overline\tau^\ast(t,z^{\frac{1}{1-\gamma}}y^{\frac{\gamma}{1-\gamma}})$: $\tau_z^\ast$ is a calendar time, whereas $\overline\tau^\ast$ is an elapsed time. Hence Proposition \ref{stopping} and equivalence of the probability measures imply}
\begin{align*}
\tau_{z_n}^\ast\longrightarrow\tau_z^\ast,\qquad\mathbb P\text{-a.s.},
\end{align*}
whenever $z_n\to z$. For $h>0$, optimality gives
\begin{align*}
\frac{\overline J(t,z+h,m,y;\tau_z^\ast)-\overline J(t,z,m,y;\tau_z^\ast)}{h}
&\leq\frac{J(t,z+h,m,y)-J(t,z,m,y)}{h}\\
&\leq\frac{\overline J(t,z+h,m,y;\tau_{z+h}^\ast)-\overline J(t,z,m,y;\tau_{z+h}^\ast)}{h}.
\end{align*}
The analogous inequalities hold for $h<0$. Using (\ref{envelopederivative}), the a.s. convergence of the stopping times, and dominated convergence, both bounds tend to
\begin{align*}
\mathbb E_{t,z,m,y}\bigg[-\int_t^{\tau_z^\ast}\xi_{s,t}\mathcal I_1^u(Z_s)ds+\xi_{\tau_z^\ast,t}\big(Q_z(Z_{\tau_z^\ast},M_{\tau_z^\ast})-g_{\tau_z^\ast}\big)\bigg].
\end{align*}
Lemma \ref{convex} ensures differentiability, so this limit equals $J_z(t,z,m,y)$ and proves the claim. In particular,
\begin{align*}
\overline{\mathcal X}(t,z^\ast,m,y)=w+g_t.
\end{align*}

Define the terminal wealth at retirement by
\begin{align*}
\phi^\ast:=-Q_z(Z_{\tau^\ast}^\ast,M_{\tau^\ast}).
\end{align*}
By (\ref{w}), $\phi^\ast>0$. Moreover, $\phi^\ast$ is $\mathcal F_{\tau^\ast}$-measurable because $Z^\ast$ and $M$ are adapted and continuous and $Q_z$ is continuous. The identity just proved yields the exact budget equality
\begin{align*}
\mathbb E_{t,w,m,y}\bigg[\int_t^{\tau^\ast}\xi_{s,t}c_s^\ast ds+\xi_{\tau^\ast,t}(\phi^\ast+g_{\tau^\ast})\bigg]=w+g_t.
\end{align*}
Lemma \ref{budegt} therefore provides an $\mathbb F$-progressively measurable portfolio process $\pi^\ast$ such that the pre-retirement wealth satisfies
\begin{align*}
W_s^\ast\geq-g_s,\qquad t\leq s<\tau^\ast,\qquad W_{\tau^\ast}^\ast=\phi^\ast.
\end{align*}

It remains to define the controls after retirement. Since $\phi^\ast=-Q_z(Z_{\tau^\ast}^\ast,M_{\tau^\ast})$, the uniqueness statement in Theorem \ref{dualityrelation} shows that the post-retirement multiplier associated with $(\phi^\ast,M_{\tau^\ast})$ is $Z_{\tau^\ast}^\ast$. Hence, conditionally on $\mathcal F_{\tau^\ast}$, we continue the process $Z^\ast$ and set, for $s\geq\tau^\ast$,
\begin{align*}
c_s^\ast:=\mathcal I_2^u(Z_s^\ast)=K^{\frac{1-\gamma}{\gamma}}(Z_s^\ast)^{-\frac{1}{\gamma}},\qquad
\pi_s^\ast:=\frac{\theta}{\sigma}Z_s^\ast Q_{zz}(Z_s^\ast,M_s).
\end{align*}
By Theorem \ref{dualityrelation}, these controls are progressively measurable and admissible after $\tau^\ast$, and the associated wealth process is $W_s^\ast=-Q_z(Z_s^\ast,M_s)\geq0$. Thus, pasting the pre-retirement and post-retirement controls gives $(c^\ast,\pi^\ast,\tau^\ast)\in\mathcal A(t,w,m,y)$.

We next show that this admissible strategy attains the dual upper bound. Multiplying the budget equality by $z^\ast$ and using
\begin{align*}
e^{-\int_t^s(\rho+M_u)du}Z_s^\ast=z^\ast\xi_{s,t},\qquad U_1(c_s^\ast)-Z_s^\ast c_s^\ast=U_1^\ast(Z_s^\ast),
\end{align*}
together with the equality case of the post-retirement duality relation,
\begin{align*}
\widehat V(\phi^\ast,M_{\tau^\ast})=Q(Z_{\tau^\ast}^\ast,M_{\tau^\ast})+Z_{\tau^\ast}^\ast\phi^\ast,
\end{align*}
we obtain
\begin{align*}
J(t,z^\ast,m,y)+z^\ast(w+g_t)=\mathbb E_{t,w,m,y}\bigg[\int_t^{\tau^\ast}e^{-\int_t^s(\rho+M_u)du}U_1(c_s^\ast)ds+e^{-\int_t^{\tau^\ast}(\rho+M_u)du}\widehat V(\phi^\ast,M_{\tau^\ast})\bigg].
\end{align*}
The right-hand side is the payoff generated by the admissible strategy constructed above. Hence,
\begin{align*}
V(t,w,m,y)\geq J(t,z^\ast,m,y)+z^\ast(w+g_t).
\end{align*}
Together with the upper bound and the fact that $z^\ast$ minimizes $J(t,z,m,y)+z(w+g_t)$, this proves
\begin{align*}
V(t,w,m,y)=\inf_{z>0}[J(t,z,m,y)+z(w+g_t)]=J(t,z^\ast,m,y)+z^\ast(w+g_t).
\end{align*}
For any fixed $z>0$, choose $w_z:=-J_z(t,z,m,y)-g_t>-g_t$. The preceding construction gives equality at $w_z$, and therefore
\begin{align*}
J(t,z,m,y)=\sup_{w>-g_t}[V(t,w,m,y)-z(w+g_t)].
\end{align*}

It remains to identify the feedback form of the optimal wealth and portfolio before $\tau^\ast$. By the conditional budget representation in Lemma \ref{budegt}, the definition of $\overline{\mathcal X}$, and the strong Markov property applied to the shifted process starting from $(s,Z_s^\ast,M_s,Y_s)$, on $\{s<\tau^\ast\}$ we have
\begin{align*}
W_s^\ast+g_s=\overline{\mathcal X}(s,Z_s^\ast,M_s,Y_s)=-J_z(s,Z_s^\ast,M_s,Y_s).
\end{align*}
Thus,
\begin{align*}
W_s^\ast=-J_z(s,Z_s^\ast,M_s,Y_s)-g_s,\qquad t\leq s<\tau^\ast.
\end{align*}
\blue{Fix the mortality path $\overline m(s)=me^{a(s-t)}$. Since $\mathcal J_x$ is smooth in the interior of the waiting region $\overline{\mathcal W}$ (cf.\ Proposition \ref{regularity}), (\ref{r1})--(\ref{Jzformula}) imply that $(s,z,y)\mapsto J_z(s,z,\overline m(s),y)$ is $C^{1,2}$ in the corresponding waiting region.  We can therefore apply It\^o's formula locally, for $t\leq s<\tau^\ast$, to
\begin{align*}
W_s^\ast=-J_z(s,Z_s^\ast,\overline m(s),Y_s)-g_s.
\end{align*}
Comparing its diffusion coefficient with (\ref{sde}) gives
\begin{align*}
\sigma\pi_s^\ast=\theta Z_s^\ast J_{zz}(s,Z_s^\ast,M_s,Y_s)-\sigma_yY_s\big(q_s+J_{zy}(s,Z_s^\ast,M_s,Y_s)\big).
\end{align*}
}
Letting the compact subsets exhaust the continuation region yields the stated portfolio rule for $t\leq s<\tau^\ast$. \blue{This identifies the already admissible portfolio constructed through Lemma \ref{budegt}.} The consumption rule holds by construction, and the post-retirement controls are those identified above. This completes the proof.
\end{proof}

\section{Recursive integration method}\label{numerical}
\renewcommand\theequation{B.\arabic{equation}}
\renewcommand\thesubsection{B.\arabic{subsection}}
The numerical method is inspired by \citet{huang1996pricing} and \citet{jeon2018portfolio}. Here we only consider the case $\gamma>1$. \blue{For the current state $(t,m)$, we have  $\overline b(t+s)=b(t+s,me^{as})=b(t+s,M_s)$. Thus (\ref{integralb2}) can be written in the original current-mortality notation as}
\begin{align}\label{6-1}
0= \int^{T-t}_0 e^{-\kappa s}\widehat{\mathbb{E}}_{\widetilde{b}(0),m}\Big[ \Big((1-K^{\frac{1-\gamma}{\gamma}}){U}^\ast_1(X_s) +1\Big) \mathds{1}_{\{X_s \leq \widetilde{b}(s)\}}\Big]ds,
\end{align}
where $\widetilde{b}(s):= b(t+s,M_s)$. Now we compute the expectation inside the integral in the right-hand side of (\ref{6-1}),
\begin{align*}
	&\widehat{\mathbb{E}}_{\widetilde{b}(0),m}\Big[ \Big((1-K^{\frac{1-\gamma}{\gamma}}){U}^\ast_1(X_s) +1\Big) \mathds{1}_{\{X_s \leq \widetilde{b}(s)\}}\Big]\\&= (1-K^{\frac{1-\gamma}{\gamma}})\frac{\gamma}{1-\gamma}\widehat{\mathbb{E}}_{\widetilde{b}(0),m}\Big[ X_s^{\frac{\gamma-1}{\gamma}}  \mathds{1}_{\{X_s \leq \widetilde{b}(s)\}}\Big]+ \widehat{\mathbb{P}}_{\widetilde{b}(0),m}[ X_s \leq \widetilde{b}(s)]\\
	&=(1-K^{\frac{1-\gamma}{\gamma}})\frac{\gamma}{1-\gamma}\widehat{\mathbb{E}}_{\widetilde{b}(0),m}\Big[ {\widetilde{b}(0)}^{\frac{\gamma-1}{\gamma}}H_sN(s,m)  \mathds{1}_{\{X_s \leq \widetilde{b}(s)\}}\Big]+ \widehat{\mathbb{P}}_{\widetilde{b}(0),m}[  X_s \leq \widetilde{b}(s)]\\	
	&=(1-K^{\frac{1-\gamma}{\gamma}})\frac{\gamma}{1-\gamma} {\widetilde{b}(0)}^{\frac{\gamma-1}{\gamma}}N(s,m)\widetilde{\mathbb{P}}_{\widetilde{b}(0),m}[ X_s \leq \widetilde{b}(s)]+ \widehat{\mathbb{P}}_{\widetilde{b}(0),m}[  X_s \leq \widetilde{b}(s)].	\end{align*}
\blue{For $s>0$, $v>0$, and $m'>0$, define
\begin{align*}
d^1(s,v,m')&:=\frac{\log v-\int^s_0[(\beta+1)(\rho-r+m'e^{au})+\mu_1-\frac{\sigma_1^2}{2}]du}{|\sigma_1|\sqrt{s}},\\
d^2(s,v,m')&:=\frac{\log v-\int^s_0[(\beta+1)(\rho-r+m'e^{au})+\mu_1+\frac{(\gamma-2)\sigma_1^2}{2\gamma}]du}{|\sigma_1|\sqrt{s}}.
\end{align*}
The factor in (\ref{N}) is
\begin{align*}
N(s,m'):=\exp\bigg(\frac{\gamma-1}{\gamma}\int^s_0\Big[(\beta+1)(\rho-r+m'e^{au})+\mu_1-\frac{\sigma_1^2}{2}\Big]du+\frac{1}{2}\Big(\frac{\gamma-1}{\gamma}\Big)^2\sigma_1^2s\bigg).
\end{align*}
Direct computation then gives
\begin{align*}
\widehat{\mathbb{P}}_{\widetilde{b}(0),m}[X_s\leq\widetilde{b}(s)]&=\Phi\bigg(d^1\Big(s,\frac{\widetilde{b}(s)}{\widetilde{b}(0)},m\Big)\bigg),\\
\widetilde{\mathbb{P}}_{\widetilde{b}(0),m}[X_s\leq\widetilde{b}(s)]&=\Phi\bigg(d^2\Big(s,\frac{\widetilde{b}(s)}{\widetilde{b}(0)},m\Big)\bigg),
\end{align*}
where $\Phi(\cdot)$ is the cumulative distribution function of a standard normal random variable. Hence, (\ref{6-1}) becomes
\begin{align*}
0=\int^{T-t}_0e^{-\kappa s}\bigg(&(1-K^{\frac{1-\gamma}{\gamma}})\frac{\gamma}{1-\gamma}{\widetilde{b}(0)}^{\frac{\gamma-1}{\gamma}}N(s,m)\Phi\bigg(d^2\Big(s,\frac{\widetilde{b}(s)}{\widetilde{b}(0)},m\Big)\bigg)\\*
&+\Phi\bigg(d^1\Big(s,\frac{\widetilde{b}(s)}{\widetilde{b}(0)},m\Big)\bigg)\bigg)ds.
\end{align*}
\blue{Set $\bar\xi:=T-t$. The terminal mortality level is $\overline m(T)=me^{a\bar\xi}$. For $u\in[0,\bar\xi]$, define}
{\color{black}
\begin{align*}
m_u&:=\overline m(T)e^{-au},\\
\widetilde b^\ast(0)&:=L,\qquad \widetilde b^\ast(u):=b(T-u,m_u),\quad 0<u\leq\bar\xi.
\end{align*}
}
Thus, $\widetilde b^\ast(\bar\xi)=b(t,m)$ and $m_{u-s}=m_ue^{as}$. For $s>0$ and $x,z,m'>0$, define
\begin{align*}
G(s,x,z,m'):=e^{-\kappa s}\bigg(&(1-K^{\frac{1-\gamma}{\gamma}})\frac{\gamma}{1-\gamma}x^{\frac{\gamma-1}{\gamma}}N(s,m')\Phi\bigg(d^2\Big(s,\frac{z}{x},m'\Big)\bigg)+\Phi\bigg(d^1\Big(s,\frac{z}{x},m'\Big)\bigg)\bigg).
\end{align*}
At $s=0$, the kernel is defined by its right limit,
\begin{align*}
G(0,x,z,m'):=\Big((1-K^{\frac{1-\gamma}{\gamma}})U_1^\ast(x)+1\Big)\bigg(\mathds{1}_{\{x<z\}}+\frac{1}{2}\mathds{1}_{\{x=z\}}\bigg).
\end{align*}
In particular, the two normal distribution functions equal $\frac{1}{2}$ at the first endpoint of the boundary recursion. For every $u\in(0,\bar\xi]$, the reversed boundary equation is
\begin{align}\label{6-2}
0=\int^u_0G(s,\widetilde b^\ast(u),\widetilde b^\ast(u-s),m_u)ds.
\end{align}
}

\blue{We divide $[0,\bar\xi]$ into $n$ equal subintervals with endpoints $\xi_j=j\Delta\xi$, where $\Delta\xi=\frac{\bar\xi}{n}$ and $j=0,1,\ldots,n$. Let $\widetilde b_j^\ast$ approximate $\widetilde b^\ast(\xi_j)$ and set \blue{$m_j:=\overline m(T)e^{-a\xi_j}$.} At the first node, the trapezoidal rule gives
\begin{align}\label{6-3}
0=\frac{\Delta\xi}{2}\big[G(0,\widetilde b_1^\ast,\widetilde b_1^\ast,m_1)+G(\xi_1,\widetilde b_1^\ast,\widetilde b_0^\ast,m_1)\big].
\end{align}
Since $\widetilde b_0^\ast=L$, only $\widetilde b_1^\ast$ is unknown. We first locate an interval in $[L,\infty)$ on which the discrete residual changes sign and then apply the bisection method; the upper endpoint is enlarged if necessary. At the second node,
\begin{align}\label{6-4}
0=\frac{\Delta\xi}{2}\big[G(0,\widetilde b_2^\ast,\widetilde b_2^\ast,m_2)+2G(\xi_1,\widetilde b_2^\ast,\widetilde b_1^\ast,m_2)+G(\xi_2,\widetilde b_2^\ast,\widetilde b_0^\ast,m_2)\big].
\end{align}
Because $\widetilde b_1^\ast$ is already known, the same procedure determines $\widetilde b_2^\ast$. More generally, for $k=2,3,\ldots,n$, $\widetilde b_k^\ast$ is obtained from
\begin{align*}
0=\frac{\Delta\xi}{2}\bigg[G(0,\widetilde b_k^\ast,\widetilde b_k^\ast,m_k)+2\sum^{k-1}_{j=1}G(\xi_j,\widetilde b_k^\ast,\widetilde b_{k-j}^\ast,m_k)+G(\xi_k,\widetilde b_k^\ast,\widetilde b_0^\ast,m_k)\bigg].
\end{align*}
\blue{Once $\{\widetilde b_j^\ast\}_{j=0}^n$ has been computed, (\ref{integralJ2}) along the same mortality trajectory yields}
\begin{align*}
\widehat{J}(t,x,m)\approx\widehat{J}_n(t,x,m):=\frac{\Delta\xi}{2}\bigg[G(0,x,\widetilde b_n^\ast,m)+2\sum^{n-1}_{j=1}G(\xi_j,x,\widetilde b_{n-j}^\ast,m)+G(\bar\xi,x,\widetilde b_0^\ast,m)\bigg].
\end{align*}
%The recursive quadrature is motivated by \citet{huang1996pricing} and \citet{jeon2018portfolio}.
}

\section*{Acknowledgments}
\blue{We thank the Editor-in-Chief, the Associate Editor, and two anonymous reviewers for their careful reading and constructive comments.}
Funded by the Deutsche Forschungsgemeinschaft (DFG, German Research Foundation) --- Project-ID 317210226 - SFB 1283. The work of Shihao Zhu was also supported by the program of China Scholarships Council.

\bibliographystyle{informs2014}
\bibliography{optimalretirement1}

@book{aksamit2017enlargement,
  title={Enlargement of filtration with finance in view},
  author={Aksamit, Anna and Jeanblanc, Monique},
   year={2017},
  publisher={Springer}
}

@article{jeon2023horizon,
  title={Horizon effect on optimal retirement decision},
  author={Jeon, Junkee and Kwak, Minsuk and Park, Kyunghyun},
  journal={Quantitative Finance},
  volume={23},
  number={1},
  pages={123--148},
  year={2023},
  publisher={Taylor \& Francis}
}

@book{krylov2008lectures,
  title={Lectures on Elliptic and Parabolic Equations in {S}obolev Spaces},
  author={Krylov, N},
  journal={Graduate Studies in Mathematics},
  volume={96},
  year={2008},
  publisher={American Mathematical Society}
}

@book{baldi2017stochastic,
  title={Stochastic calculus},
  author={Baldi, Paolo},
  year={2017},
  publisher={Springer}
}

@article{bensoussan2016unemployment,
  title={Unemployment risks and optimal retirement in an incomplete market},
  author={Bensoussan, Alain and Jang, Bong-Gyu and Park, Seyoung},
  journal={Operations Research},
  volume={64},
  number={4},
  pages={1015--1032},
  year={2016},
  publisher={INFORMS}
}

@article{he2022optimal,
  title={Optimal asset allocation, consumption and retirement time with the variation in habitual persistence},
  author={He, Lin and Liang, Zongxia and Song, Yilun and Ye, Qi},
  journal={Insurance: Mathematics and Economics},
  volume={102},
  pages={188--202},
  year={2022},
  publisher={Elsevier}
}

@article{banks1998there,
  title={Is there a retirement-savings puzzle?},
  author={Banks, James and Blundell, Richard and Tanner, Sarah},
  journal={American Economic Review},
  number = {4},
  pages={769--788},
  volume = {88},
  year={1998},
  publisher={JSTOR}
}

@article{jeon2018portfolio,
  title={Portfolio selection with consumption ratcheting},
  author={Jeon, Junkee and Koo, Hyeng Keun and Shin, Yong Hyun},
  journal={Journal of Economic Dynamics and Control},
  volume={92},
  pages={153--182},
  year={2018},
  publisher={Elsevier}
}

@article{chen2021retirement,
  title={On retirement time decision making},
  author={Chen, An and Hentschel, Felix and Steffensen, Mogens},
  journal={Insurance: Mathematics and Economics},
  volume={100},
  pages={107--129},
  year={2021},
  publisher={Elsevier}
}

@book{peskir2006optimal,
  title={Optimal stopping and free-boundary problems},
  author={Peskir, Goran and Shiryaev, Albert},
  year={2006},
  publisher={Springer}
}

@article{jin2006disutility,
  title={Disutility, optimal retirement, and portfolio selection},
  author={Choi, Kyoung and Shim, Gyoocheol},
  journal={Mathematical Finance},
  volume={16},
  number={2},
  pages={443--467},
  year={2006},
  publisher={Wiley Online Library}
}

@article{yang2018optimal,
  title={Optimal consumption and portfolio selection with early retirement option},
  author={Yang, Zhou and Koo, Hyeng Keun},
  journal={Mathematics of Operations Research},
  volume={43},
  number={4},
  pages={1378--1404},
  year={2018},
  publisher={INFORMS}
}

@article{yang2021optimal,
  title={Optimal retirement in a general market environment},
  author={Yang, Zhou and Koo, Hyeng Keun and Shin, Yong Hyun},
  journal={Applied Mathematics \& Optimization},
  volume={84},
  number={1},
  pages={1083--1130},
  year={2021},
  publisher={Springer}
}

@article{chen2021optimal,
  title={Optimal retirement under partial information},
  author={Chen, Kexin and Jeon, Junkee and Wong, Hoi Ying},
  journal={Mathematics of Operations Research},
  volume={47},
  number={3},
  pages={1802--1832},
  year={2022},
  publisher={INFORMS}
}

@article{jeon2020optimal,
  title={Optimal retirement and portfolio selection with consumption ratcheting},
  author={Jeon, Junkee and Park, Kyunghyun},
  journal={Mathematics and Financial Economics},
  volume={14},
  number={3},
  pages={353--397},
  year={2020},
  publisher={Springer}
}

@article{bodie2004optimal,
  title={Optimal consumption--portfolio choices and retirement planning},
  author={Bodie, Zvi and Detemple, J{\'e}r{\^o}me B and Otruba, Susanne and Walter, Stephan},
  journal={Journal of Economic Dynamics and Control},
  volume={28},
  number={6},
  pages={1115--1148},
  year={2004},
  publisher={Elsevier}
}

@article{dybvig2010lifetime,
  title={Lifetime consumption and investment: {r}etirement and constrained borrowing},
  author={Dybvig, Philip H and Liu, Hong},
  journal={Journal of Economic Theory},
  volume={145},
  number={3},
  pages={885--907},
  year={2010},
  publisher={Elsevier}
}

@article{guan2025retirement,
  title={Retirement decision with addictive habit persistence in a jump diffusion market},
  author={Guan, Guohui and Huang, Qitao and Liang, Zongxia and Yuan, Fengyi},
  journal={SIAM Journal on Financial Mathematics},
  volume={16},
  number={3},
  pages={912--958},
  year={2025},
  publisher={SIAM}
}

@book{karatzas1998methods,
  title={Methods of mathematical finance},
  author={Karatzas, Ioannis and Shreve, Steven E},
  volume={39},
  year={1998},
  publisher={Springer}
}

@article{karatzas2000utility,
  title={Utility maximization with discretionary stopping},
  author={Karatzas, Ioannis and Wang, Hui},
  journal={SIAM Journal on Control and Optimization},
  volume={39},
  number={1},
  pages={306--329},
  year={2000},
  publisher={SIAM}
}

@article{farhi2007saving,
  title={Saving and investing for early retirement: A theoretical analysis},
  author={Farhi, Emmanuel and Panageas, Stavros},
  journal={Journal of Financial Economics},
  volume={83},
  number={1},
  pages={87--121},
  year={2007},
  publisher={Elsevier}
}

@article{choi2008optimal,
  title={Optimal portfolio, consumption-leisure and retirement choice problem with {CES} utility},
  author={Choi, Kyoung Jin and Shim, Gyoocheol and Shin, Yong Hyun},
  journal={Mathematical Finance},
  volume={18},
  number={3},
  pages={445--472},
  year={2008},
  publisher={Wiley Online Library}
}

@article{chen2018optimal,
  title={Optimal retirement time under habit persistence: {W}hat makes individuals retire early?},
  author={Chen, An and Hentschel, Felix and Xu, Xian},
  journal={Scandinavian Actuarial Journal},
  volume={2018},
  number={3},
  pages={225--249},
  year={2018},
  publisher={Taylor \& Francis}
}

@book{karatzas2012brownian,
  title={Brownian motion and stochastic calculus},
  author={Karatzas, Ioannis and Shreve, Steven E},
  volume={113},
  year={1998},
  publisher={Springer Science \& Business Media}
}

@article{dybvig2011verification,
  title={Verification theorems for models of optimal consumption and investment with retirement and constrained borrowing},
  author={Dybvig, Philip H and Liu, Hong},
  journal={Mathematics of Operations Research},
  volume={36},
  number={4},
  pages={620--635},
  year={2011},
  publisher={INFORMS}
}

@article{de2017dividend,
  title={The dividend problem with a finite horizon},
  author={De Angelis, Tiziano and Ekstr{\"o}m, Erik},
  journal={The Annals of Applied Probability},
  volume={27},
  number={6},
  pages={3525--3546},
  year={2017},
  publisher={Institute of Mathematical Statistics}
}

@book{bensoussan2011applications,
  title={Applications of variational inequalities in stochastic control},
  author={Bensoussan, Alain and Lions, J-L},
  year={1982},
  publisher={NorthHolland
(Amsterdam)}
}

@article{peskir2005american,
  title={On the {American} option problem},
  author={Peskir, Goran},
  journal={Mathematical Finance},
  volume={15},
  number={1},
  pages={169--181},
  year={2005},
  publisher={Wiley Online Library}
}

@article{gompertz1825nature,
  title={On the nature of the function expressive of the law of human mortality, and on a new mode of determining the value of life contingencies},
  author={Gompertz, B},
  journal={Philos. Trans. Roy. Soc. London},
  volume={115},
  pages={513--585},
  year={1825}
}

@article{de2019lipschitz,
  title={On {Lipschitz} continuous optimal stopping boundaries},
  author={De Angelis, Tiziano and Stabile, Gabriele},
  journal={SIAM Journal on Control and Optimization},
  volume={57},
  number={1},
  pages={402--436},
  year={2019},
  publisher={SIAM}
}

@article{szinovacz2014recession,
  title={Recession and expected retirement age: Another look at the evidence},
  author={Szinovacz, Maximiliane E and Martin, Lauren and Davey, Adam},
  journal={The Gerontologist},
  volume={54},
  number={2},
  pages={245--257},
  year={2014},
  publisher={Oxford University Press US}
}

@article{wang2014psychological,
  title={Psychological research on retirement},
  author={Wang, Mo and Shi, Junqi},
  journal={Annual Review of Psychology},
  volume={65},
  pages={209--233},
  year={2014},
  publisher={Annual Reviews}
}

@article{fisher2016retirement,
  title={Retirement timing: A review and recommendations for future research},
  author={Fisher, Gwenith G and Chaffee, Dorey S and Sonnega, Amanda},
  journal={Work, Aging and Retirement},
  volume={2},
  number={2},
  pages={230--261},
  year={2016},
  publisher={Oxford University Press}
}

@article{de2019free,
  title={On the free boundary of an annuity purchase},
  author={De Angelis, Tiziano and Stabile, Gabriele},
  journal={Finance and Stochastics},
  volume={23},
  number={1},
  pages={97--137},
  year={2019},
  publisher={Springer}
}

@article{milevsky2007annuitization,
  title={Annuitization and asset allocation},
  author={Milevsky, Moshe A and Young, Virginia R},
  journal={Journal of Economic Dynamics and Control},
  volume={31},
  number={9},
  pages={3138--3177},
  year={2007},
  publisher={Elsevier}
}

@article{honig1998married,
  title={Married women's retirement expectations: Do pensions and social security matter?},
  author={Honig, Marjorie},
  journal={American Economic Review},
  volume={88},
  number={2},
  pages={202--206},
  year={1998},
  publisher={JSTOR}
}

@article{jang2024optimal,
  title={Optimal investment, heterogeneous consumption, and best time for retirement},
  author={Jang, Hyun Jin and Xu, Zuo Quan and Zheng, Harry},
  journal={Operations Research},
  volume={72},
  number={2},
  pages={832--847},
  year={2024},
  publisher={INFORMS}
}

@article{park2023robust,
  title={Robust Retirement with Return Ambiguity: Optimal {G}-Stopping Time in Dual Space},
  author={Park, Kyunghyun and Wong, Hoi Ying},
  journal={SIAM Journal on Control and Optimization},
  volume={61},
  number={3},
  pages={1009--1037},
  year={2023},
  publisher={SIAM}
}

@article{huang1996pricing,
  title={Pricing and hedging {A}merican options: a recursive integration method},
  author={Huang, Jing-zhi and Subrahmanyam, Marti G and Yu, G George},
  journal={The Review of Financial Studies},
  volume={9},
  number={1},
  pages={277--300},
  year={1996},
  publisher={Oxford University Press}
}

\end{document}